\documentclass{article}
\usepackage{tikz}
\usepackage{amsmath}
\usepackage{amsfonts}
\usepackage{amssymb}
\usepackage{xspace}
\usepackage{graphicx}
\numberwithin{equation}{section}
\makeatletter
\newenvironment{prooof}{\par\noindent{\em Proof:}
}{\hfill\llap{$\Box$}\vspace{1\baselineskip}\par\noindent}

\newenvironment{proofof}{\par\noindent{\em Proof}
}{\hfill\llap{$\Box$}\vspace{1\baselineskip}\par\noindent}

\newcommand{\adl}{\vspace{1\baselineskip}}

\newtheorem{theorem}{Theorem}[section]
\newtheorem{proposition}[theorem]{Proposition}
\newtheorem{lemma}[theorem]{Lemma}
\newtheorem{corollary}[theorem]{Corollary}
\newtheorem{remark}[theorem]{Remark}
\newtheorem{definition}[theorem]{Definition}
\newtheorem{example}[theorem]{Example}

\newcommand{\beq}{\begin{equation}}
\newcommand{\eeq}{\end{equation}}

\newcommand{\ba}{\begin{array}}
\newcommand{\ea}{\end{array}}
\newcommand{\bt}{\begin{theorem}}
\newcommand{\et}{\end{theorem}}
\newcommand{\bp}{\begin{proposition}}
\newcommand{\ep}{\end{proposition}}
\newcommand{\bl}{\begin{lemma}}
\newcommand{\el}{\end{lemma}}
\newcommand{\bc}{\begin{corollary}}
\newcommand{\ec}{\end{corollary}}
\newcommand{\bi}{\begin{itemize}}
\newcommand{\ei}{\end{itemize}}
\newcommand{\ben}{\begin{enumerate}}
\newcommand{\een}{\end{enumerate}}
\newcommand{\bpf}{\begin{prooof}}
\newcommand{\epf}{\end{prooof}}
\newcommand{\bpff}{\begin{proofof}}
\newcommand{\epff}{\end{proofof}}
\newcommand{\bdf}{\begin{definition}\rm}
\newcommand{\edf}{\end{definition}}
\newcommand{\br}{\begin{remark}\rm}
\newcommand{\er}{\end{remark}}
\newcommand{\bex}{\begin{example}\rm}
\newcommand{\eex}{\end{example}}
\def\pri{\hbox to 10pt{\hfil\hbox to 0.4pt{\vrule height5pt width0.4pt
                 depth0pt}\vrule width5pt height0.4pt depth0pt\hfil}}
\newcommand{\TC}{{\rm TC}}
\newcommand{\CT}{{\rm CT}}

\newcommand{\TAT}{{\rm TAT}}

\newcommand{\gk}{{\bf{k}}}

\newcommand{\gK}{{\bf K}}
\newcommand{\gN}{{\bf N}}

\newcommand{\gT}{{\bf T}}

\newcommand{\calL}{{\cal L}}

\newcommand{\gR}{{\mathbb R}}

\newcommand{\Nat}{{\mathbb N}}

\newcommand{\F}{{\cal F}}

\newcommand{\BV}{\mathop{\rm BV}\nolimits}

\def\sgn{\mathop{\rm sgn}\nolimits}

\newcommand{\wc}{\rightharpoonup}

\newcommand{\RP}{{\mathbb{RP}}^2}

\newcommand{\SP}{{{\mathbb S}^2}}

\newcommand{\SN}{{{\mathbb S}^{N}}}

\newcommand{\Sjp}{{{\mathbb S}^{j}_{p}}}
\newcommand{\RPN}{{{\mathbb{RP}}^{N}}}

\newcommand{\RPjp}{{{\mathbb{RP}}^{j}_{p}}}

\def\mesh{\mathop{\rm mesh}\nolimits}
\def\span{\mathop{\rm span}\nolimits}

\newcommand{\gtau}{{\boldsymbol{\tau}}}

\newcommand{\bfg}{{\boldsymbol{\gamma}}}

\newcommand{\gt}{{\bf t}}
\newcommand{\gn}{{\bf n}}
\newcommand{\gb}{{\bf b}}
\newcommand{\gu}{{\bf u}}
\newcommand{\gv}{{\bf v}}
\newcommand{\go}{{\bf o}}

\newcommand{\gc}{{\bf c}}
\newcommand{\ga}{{\bf a}}

\newcommand{\Var}{\mathop{\rm Var}\nolimits}

\let\a=\alpha
\let\be=\beta

\let\d=\delta
\let\e=\varepsilon

\let\vf=\varphi
\let\g=\gamma

\let\l=\lambda
\let\m=\mu

\let\p=\pi

\let\s=\sigma
\let\t=\theta
\let\tt=\tau
\let\x=\xi
\let\y=\eta

\let\vf=\varphi

\let\SS=\Sigma
\let\O=\Omega

\let\wih=\widehat
\let\wdg=\wedge

\let\wid=\widetilde

\let\sb=\subset

\let\fa=\forall
\let\tim=\times

\let\sm=\setminus
\let\ol=\overline

\let\ds=\displaystyle

\let\lan=\langle
\let\ran=\rangle

\let\i=\infty
\let\lm=\limits

\title{\Large \bf Weak curvatures of irregular curves \\ in high dimension Euclidean spaces}
\author{\it Domenico Mucci and Alberto Saracco
\footnote{%
{\sc Dipartimento di Scienze Matematiche,
Fisiche ed Informatiche, Universit\`{a} di Parma,
Parco Area delle Scienze 53/A, I-43124 Parma, Italy.
E-mail: domenico.mucci@unipr.it, alberto.saracco@unipr.it}
}
}
\begin{document}
\date{}
\topskip=1.5truecm \maketitle \topskip=1.5truecm \maketitle

\maketitle
%
{\small{\bf Abstract.}
We deal with a robust notion of weak normals for a wide class of irregular curves defined in Euclidean spaces of high dimension.
Concerning polygonal curves, the discrete normals are built up through a Gram-Schmidt procedure applied to consecutive oriented segments,
and they naturally live in the projective space associated to the Gauss hyper-sphere. By using sequences of inscribed polygonals with infinitesimal modulus, a relaxed notion of total variation of the $j$-th normal to a generic curve is then introduced. For smooth curves satisfying the Jordan system, in fact, our relaxed notion agrees with the length of the smooth $j$-th normal.
Correspondingly, a good notion of weak $j$-th normal of irregular curves with finite relaxed energy is introduced, and it turns out to be the strong limit of any sequence of approximating polygonals.
The length of our weak normal agrees with the corresponding relaxed energy, for which a related integral-geometric formula is also obtained.
We then discuss a wider class of smooth curves for which the weak normal is strictly related to the classical one, outside the inflection points.
Finally, starting from the first variation of the length of the weak $j$-th normal, a natural notion of curvature measure is also analyzed.}
\adl\par\noindent
{\bf Keywords :} Jordan system; Relaxed energies; Polygonals; Non-smooth curves
\adl\par\noindent
{\bf MSC :} {53A04; 49J45}
\section{Introduction}
\label{intro}
The well-known notions of curvature and torsion of a smooth rectifiable curve $\gc$ in $\gR^3$ were independently obtained by Frenet and Serret.
%
%
The extension to smooth curves in high dimension Euclidean spaces $\gR^{N+1}$, where $N\geq 3$, goes back to the contribution by C.~Jordan \cite{J},
who noticed that by applying the Gram-Schmidt procedure to the independent vectors $\dot\gc(s),\,\gc^{(2)}(s),\ldots,\gc^{(N)}(s)$ one obtains a moving frame ${\bf e}(s):=(\gt(s),\gn_1(s),\ldots,\gn_N(s))$ along the curve, where $\gt$ is the tantrix (or tangent indicatrix) and $\gn_j$ is the $j$-th curvature, for $j=1,\ldots,N$.
Assuming $\gc$ parameterized by arc-length $s$, the Jordan system $\dot{\bf e}(s)=F(s)\,{\bf e}(s)$ involves a skew-symmetric and tri-diagonal square matrix $F(s)$ of order $N+1$, whose entries depend on the curvature functions $\gk_j(s)$, where $j=1,\ldots,N$.
\par In this framework, H.~Gluck \cite{Gl} produced an algorithm for computing the higher order curvatures, whereas more recently E.~Gutkin \cite{Gu} studied curvature estimates, natural invariants, and discussed the case of curves contained in Riemannian manifolds and homogeneous spaces.
\smallskip\par In this paper, we are interested in analyzing an analogous theory concerning irregular curves.
The main historical contribution goes back to the work by A.~D.~Alexandrov and Yu.~G.~Reshetnyak \cite{AR}. In the last section of his more recent survey paper \cite{Re}, Reshetnyak also discussed possible ways to extend their theory of irregular curves to the high codimension case.
\par To this purpose, we recall that the definition of {\em complete torsion} $\CT(P)$ of polygonals $P$ in $\gR^3$ given by Alexandrov-Reshetnyak \cite{AR},
who essentially take the distance in $\SP$ between consecutive discrete binormals, implies that planar polygonals may have positive torsion at ``inflections points".
Defining the complete torsion $\CT(\gc)$ of curves $\gc$ in $\gR^3$ as the supremum of the complete torsion among the inscribed polygonals,
they obtain in \cite[p.~244]{AR} that any curve with finite complete torsion and with no points of return must have finite {\em total curvature} $\TC(\gc)$, see \eqref{TC}.
\par Notice however that a rectifiable smooth curve in $\gR^3$ may have unbounded total curvature but zero torsion (just consider a planar curve).
On the other hand, the (absolute value of the) torsion may be seen as the curvature of the tantrix, when computed in the sense of spherical geometry.
\par For these reasons, in our paper \cite{MStor} on irregular curves in $\gR^3$, following the approach by M.~A.~Penna \cite{Pe}, we
defined the binormal indicatrix $\gb_P$ of a polygonal $P$ as
the arc-length parameterization of the polar in the projective plane $\RP$ of the tantrix $\gt_P$.
Therefore, the {\em total absolute torsion} $\TAT(P)$ of $P$ is equal to the length of the curve $\gb_P$ in $\RP$.
Furthermore, by exploiting the polarity in $\RP$, we also discussed a notion of principal normal $\gn_P$.
We remark that a similar definition has been introduced by T.~F.~Banchoff in his paper \cite{Ba} on space polygons.
\bigskip\par\noindent{\large\sc Content of the paper.}
When dealing with polygonal curves $P$ in high dimension Euclidean spaces, the polarity argument previously described fails to hold. Therefore,
in this paper we follow a different approach, based on the orthonormalization procedure.
Referring to Sec.~\ref{Sec:discr} for details on the construction,
in order to define the discrete $j$-th normal to a polygonal $P$, for $j\leq N-1$, we consider lists of $j+1$ consecutive segments of $P$ that do not lay on any affine $j$-space of $\gR^{N+1}$.
Therefore, they define a {\em discrete osculating $(j+1)$-space}, and we choose the last unit vector in $\SN$ obtained by means of the Gram-Schmidt procedure.
We then consider the corresponding points in the projective space $\RPN$, that are naturally ordered w.r.t. the consecutive segments of the polygonal $P$, and
define the {\em $j$-th normal $[\gn_j](P)$} as the curve in $\RPN$ obtained by connecting these consecutive points with geodesic arcs.
\par As to the last normal $[\gn_N](P)$, we consider the equivalence classes in $\RPN$ of the orthogonal directions to the
discrete osculating $N$-spaces, and argue the same way as above. Therefore, when $N=2$, we recover our notion of binormal indicatrix $\gb_P$ from \cite{MStor}.
\par In Theorem~\ref{Tappr}, we show that for any {\em smoothly turning curve} $\gc$ we
can find a sequence $\{P_n\}$ of inscribed polygonals, with $\mesh P_n\to 0$, such that the length $\calL_\RPN([\gn_j](P_n))$ of the discrete $j$-th normal to $P_n$ converges to the length $\calL_\SN(\gn_j)$ of the $j$-th normal $\gn_j$ to the curve $\gc$, i.e.,
$$ \lim_{n\to\i}\calL_\RPN([\gn_j](P_n))=\int_a^b\Vert\dot\gn_j(s)\Vert\,ds \qquad \fa\,j=1,\ldots,N\,. $$
\par We recall that by the Jordan formulas~\eqref{FSN}, one has $ \Vert\dot\gn_j(s)\Vert=\sqrt{\gk_{j}^2(s)+\gk_{j+1}^2(s) }$, if $j<N$, whereas
$ \Vert\dot\gn_N(s)\Vert=|\gk_N(s)|$ when $j=N$, for the last normal.
\par A smoothly turning curve $\gc$, see Definition~\ref{Dst}, essentially corresponds to the
regular curves considered by Gutkin \cite{Gu}, and it satisfies the Jordan system~\eqref{FSN}.
In order to construct the approximating sequence $\{P_n\}$, in Sec.~\ref{Sec:GS}, at a given interior point $\gc(s)$ we consider the inscribed polygonals corresponding to vertexes at arc-length distance $h>0$, see \eqref{vjh}. The $j$-th normals of such polygonals can be written in terms of the Taylor expansions of $\gc(s)$, see Propositions~\ref{PGMN} and \ref{PGM4},
where computations are postponed to the appendix.
\smallskip\par Motivated by the previous density result, in Sec.~\ref{Sec:rel},
we introduce a {\em relaxed notion of total variation of the $j$-th normal} to a generic curve $\gc$ in $\gR^{N+1}$.
Now, differently to what happens for length and total curvature, the monotonicity formula fails to hold in general for the length of the discrete
$j$-th normal to polygonals, see Remark~\ref{Rmon}, Example~\ref{Emon}, and Figure~\ref{fig:controex}.
Therefore, we are led to follow the approach introduced by Alexandrov-Reshetnyak \cite{AR}, that involves the notion of {\em modulus} $\m_\gc(P)$ of a polygonal $P$ inscribed in $\gc$, say $P\prec \gc$, and we define:
$$ \F_j(\gc):=\lim_{\e\to 0^+}\sup\{\calL_\RPN([\gn_j](P)) \mid P \prec \gc\,,\,\,\m_c(P)<\e\} \qquad j=1,\ldots,N \,. $$
\par We point out, in fact, that for polygonal curves $P$ in $\gR^{N+1}$ one has
$$
 \F_j(P) = \calL_\RPN([\gn_j](P)) \qquad\fa\,j=1,\ldots,N $$
whereas in the case $N=2$, the relaxed total variation of the last normal agrees with the total absolute torsion of curves $\gc$ in $\gR^3$ that we analyzed in \cite{MStor}.
\par Most importantly, in Proposition~\ref{PFjc} we show that if a curve $\gc$ satisfies $\F_{j}(\gc)<\i$ and $\F_{j-1}(\gc)<\i$ for some $j=2,\ldots,N$,
then for any sequence $\{P_n\}$ of inscribed polygonals for which $\m_\gc(P_n)\to 0$ one has:
$$ \lim_{n\to\i} \calL_\RPN([\gn_j](P_n))=\F_j(\gc)\,. $$
In the case $j=1$, the same conclusion holds true for any curve $\gc$ satisfying $\TC(\gc)<\i$.
\par Therefore, for smoothly turning curves, in Proposition~\ref{PFjcsmooth} we also obtain the explicit formulas:
$$ \F_j(\gc)=\int_a^b\Vert\dot\gn_j(s)\Vert\,ds\,. $$
%
\noindent{\large\sc Weak normals.} The previous continuity property is a consequence of the Main Result of this paper, Theorem~\ref{Tjcurv}, that
justifies our notion of {\em weak $j$-th normal} $[\gn_j](\gc)$ to a curve $\gc$.
Notice that in this paper we do not need to restrict to consider simple curves, since our construction is based on local arguments.
\smallskip\par More precisely, we have:
\smallskip\par\noindent
{\bf Main Result.} {\em Let $N\geq 2$ and $\gc$ be a curve in $\gR^{N+1}$ such that $\F_{j}(\gc)<\i$ and $\F_{j-1}(\gc)<\i$ for some $j=2,\ldots,N$.
There exists a rectifiable curve $[\gn_j](\gc):[0,L_j]\to\RPN$ parameterized by arc-length, where
$$L_j:=\F_{j}(\gc)= \calL_{\RPN}([\gn_j](\gc)) $$
satisfying the following property.
For any sequence $\{P_n\}$ of inscribed polygonal curves, let $\bfg^j_n:[0,L_j]\to\RPN$ denote for each $n$ the parameterization with constant velocity
of the discrete $j$-th normal  $[\gn_j](P_n)$ to $P_n$, see Definition~$\ref{Dnj}$.
If $\m_\gc(P_n)\to 0$, then $\bfg^j_n\to [\gn_j](\gc)$ uniformly on $[0,L_j]$ and
$$\calL_{\RPN}(\bfg^j_n)=\calL_{\RPN}([\gn_j](P_n))\to\calL_{\RPN}([\gn_j](\gc))$$
as $n\to\i$, where, we recall, $\calL_{\RPN}([\gn_j](P_n))=\F_j(P_n)$.
Moreover, the arc-length derivative of the curve $[\gn_j](\gc)$ is a function of bounded variation.
Finally, in the case $j=1$, for any curve
$\gc$ in $\gR^{N+1}$ satisfying $\TC(\gc)<\i$, one has $\F_1(\gc)<\i$ and the same conclusion as above holds true.}
\smallskip\par In Sec.~\ref{Sec:weakn}, the proof of our Main Result proceeds by steps. Firstly, we obtain the curve $[\gn_j](\gc)$ by means of an optimal approximating sequence,
where we have to apply the sequential weak-* compactness theorem for one-dimensional $\BV$-functions, see \cite{AFP}. We thus need a uniform bound for the
total variation of the tantrix associated to a continuous lifting of the curve $\bfg^j_n$. It holds true provided that we assume that $\F_{j-1}(\gc)<\i$, when $j>1$, in addition to the natural hypothesis $\F_{j}(\gc)<\i$, see Remark~\ref{Rcomp}.
\par Following some ideas taken from our paper \cite{MStor}, we then deal with the case $j=N$ by exploiting the polarity of the last normal.
\par In order to analyze the case $1<j<N$ of the intermediate normals, we then make use of an {\em integral-geometric formula for polygonals}, see \eqref{igP}.
It is obtained in Sec.~\ref{Sec:est}, as a consequence of our Theorem~\ref{TAR}, where we extend the integral-geometric formula for polygonal curves in $\SN$ due to Alexandrov-Reshetnyak \cite[Thm. 6.2.2, p. 190]{AR}, who only treated the case of projections onto low dimension spaces.
\par In Proposition~\ref{Pineqj}, we also obtain the following inequality concerning the total curvature of the discrete $j$-th normal to a polygonal curve $P$:
$$ \TC([\gn_j](P))\leq \calL_{\RPN}([\gn_{j-1}](P))+\calL_{\RPN}([\gn_{j}](P))\qquad \fa\,j=2,\ldots,N $$
that is crucial in the previously cited compactness argument.
\par At the final step, we treat the case of the first normal, using this time that
$$\TC([\gn_1](P))\leq\calL_{\SN}(\gt_P)+\calL_{\RPN}([\gn_{1}](P))\,,\qquad \calL_{\SN}(\gt_P)=\TC(P) $$
and that we always have $ \calL_{\RPN}([\gn_1](P))\leq\TC(P)$, see Proposition~\ref{PigP}.
\par As a consequence, if a curve $\gc$ satisfies for some integer $1\leq j\leq N-1$ the hypotheses of our main result, in Corollary~\ref{Cigc} we also obtain the
following integral-geometric formula:
$$ \F_j(\gc)=\int_{G_{j+1}\gR^{N+1}}\F_j(\p_p(\gc)))\,d\m_{j+1}(p)\,. $$
Here, $G_{j+1}\gR^{N+1}$ is the Grassmannian of
the unoriented ${(j+1)}$-planes in $\gR^{N+1}$, $\m_{j+1}$ is the corresponding Haar measure, and $\p_p$ is the
orthogonal projection of $\gR^{N+1}$ onto an element $p$ in $G_{j+1}\gR^{N+1}$.
\bigskip\par\noindent{\large\sc Other results.} In Sec.~\ref{Sec:relation}, we analyze the relationship between our weak $j$-th normal $[\gn_j](\gc)$ and the classical $j$-th normal $\gn_j$.
For smoothly turning curves, the expected result is obtained in Proposition~\ref{Csmooth}.
\par With the aim of finding a wider class of smooth curves $\gc$ satisfying a similar relationship, we point out that the main property we need to preserve is the
{\em existence and continuity of the osculating $(j+1)$-spaces}.
Such a property is guaranteed for {\em mildly smoothly turning curves} as in our Definition~\ref{Dms}, see Proposition~\ref{Posc}.
Any such curve satisfies the Jordan system~\eqref{FSN} outside a finite set $\SS$ of points, Proposition~\ref{Pcomp}.
Also, both the convergence result in Theorem~\ref{Tappr} and the representation formula in Proposition~\ref{PFjcsmooth} for the relaxed total variation of the $j$-th normal continue to hold, see Propositions~\ref{Pappr} and \ref{PFjcmild}.
Finally, the relationship between the weak $j$-th normal $[\gn_j](\gc)$ from Theorem~\ref{Tjcurv} and the smooth $j$-th normal is analyzed in Proposition~\ref{Cmild}.
\par In Sec.~\ref{Sec:curvmeas}, we deal with the measure given by the distributional derivative of the arc-length derivative of the weak $j$-th normal to a curve satisfying the hypotheses of our Main Result. The case of the tangent indicatrix was firstly discussed in \cite{CFKSW}, see also \cite{Su_curv},
where the authors introduced the notion of {\em curvature force}. It comes into the play when considering the first variation of the length of curves with finite total curvature.
When $N=2$, the {\em torsion force} was similarly discussed in \cite{MStor}, where we considered tangential variations of the length of the tantrix.
\par Roughly speaking, a continuous lifting $\bfg^j:[0,L_j]\to\SN$ of the curve $[\gn_j](\gc)$ in our Main Result is such that its arc-length derivative $\dot\bfg^j$ is a function of bounded variation, and its distributional derivative appears when computing the first variation
of the length $\calL_{\SN}(\bfg^j)$, see formula~\eqref{var}.
In particular, for smoothly turning curves, we obtain the formula
$$ D\dot\bfg^j=\vf_{j\,\#}\m_j\,,\qquad \m_j:={d\over ds}\,{\dot\gn_j(s)\over \Vert\dot\gn_j(s)\Vert}\,\calL^1\pri]a,b[ $$
where on the left-hand side we are denoting the push forward of the measure $\m_j$ by the transition function $t=\vf_j(s):=\int_a^s\Vert\dot\gn_j(\l)\Vert\,d\l$, see also Example~\ref{Emeas}.
\par Finally, the curvature measures associated to our mildly smoothly turning curves are also analyzed, yielding to more general properties.
\section{Gram-Schmidt procedure}\label{Sec:GS}
In this section, we deal with Taylor expansions of inscribed polygonals to smooth curves.
By means of a Gram-Schmidt procedure, we analyze the relationship between the approximate frame and the Jordan frame of the given curve.
For that reason, we introduce a suitable notion of {\em smoothly turning curve}, see Definition~\ref{Dst}.
We first discuss the first two normals, and then consider the general case.
\adl\par\noindent{\large\sc The first two normals.} Let $N\geq 2$ and $\gc:[a,b]\to\gR^{N+1}$ be a curve of class $C^3$ parameterized by arc-length, so that $\Vert\dot \gc\Vert=1$.
Denoting by $\gc^{(k)}$ the $k$-th arc-length derivative of $\gc$, assume that the triplet
$(\dot\gc(s),\gc^{(2)}(s),\gc^{(3)}(s))$ is linearly independent for each $s$.
%
The first two Frenet-Serret formulas give
$$ \dot\gt=\gk_1\,\gn_1\,,\quad \dot\gn_1=-\gk_1\,\gt+\gk_2\,\gn_2  $$
where $\gt:=\dot\gc\in\SN$ is the unit tangent vector, $\gk_1:=\Vert\gc^{(2)}\Vert$ is the first curvature, $\gn_1:=\gc^{(2)}/\Vert\gc^{(2)}\Vert\in\SN$ is the
first unit normal, $\gk_2\in\gR$ is the second curvature and $\gn_2\in\SN$ is the second unit normal. Notice that when $N=2$ one has $\gk_2=\gtau$, the torsion of the curve, and $\gn_2=\gb$, the binormal vector $\gb:=\gt\tim\gn$.
\par Denoting by $\bullet$ the scalar product in $\gR^{N+1}$, and following an argument that goes back to Jordan \cite{J},
we thus compute
$$ \gk_2\,\gn_2 =\gk_1\,\gt+\dot\gn_1= \Vert\gc^{(2)}\Vert\,\dot\gc+{d\over ds}\Bigl(\frac{\gc^{(2)}}{\Vert\gc^{(2)}\Vert} \Bigr)=
\frac{1}{\Vert\gc^{(2)}\Vert}\,\Bigl( \Vert\gc^{(2)}\Vert^2\,\dot\gc+\gc^{(3)}-\frac{\gc^{(2)}\bullet\gc^{(3)}}{\Vert\gc^{(2)}\Vert^2}\,\gc^{(2)}\Bigr)\,. $$
We now recall that $\dot\gc\bullet\gc^{(2)}=0$ and that $\dot\gc\bullet\gc^{(3)}=-\Vert\gc^{(2)}\Vert^2$.
Therefore, according to the Gram-Schmidt procedure one has:
$$ \gn_2 = \frac{\gc^{(3)\perp}}{\Vert\gc^{(3)\perp}\Vert}\,,\qquad
\gc^{(3)\perp}:= \gc^{(3)} - \frac{\gc^{(3)}\bullet\dot\gc}{\Vert\dot\gc\Vert^2}\,\dot\gc-\frac{\gc^{(3)}\bullet\gc^{(2)}}{\Vert\gc^{(2)}\Vert^2}\,\gc^{(2)}\,.$$
\par We now choose some point $s\in]a,b[$ and for each $h>0$ small enough we consider the three vectors
\beq\label{vh}
\ds \gv_0(h):=\frac {\gc(s+h)-\gc(s-h)}{2h}\,, \   \gv_1(h):=\frac {\gc(s-3h)-\gc(s-h)}{2h}\,,\  \gv_2(h):=\frac {\gc(s+3h)-\gc(s+h)}{2h}\,. \eeq
In the sequel, we omit to write the dependence on $s$, and denote by $\go(h^n)$ a continuous vector function such that $\Vert\go(h^n)\Vert=o(h^n)$, for each $n\in\Nat$, i.e., $\Vert\go(h^n)\Vert/h^n\to 0$ as $h\to 0$.
\par By taking the third order expansions of $\gc(s)$ and by applying the Gram-Schmidt procedure, we obtain:
\bp\label{PGM3} We have:
$$ \gt(h):=\frac{\gv_0(h)}{\Vert\gv_0(h)\Vert}=\dot \gc+\frac 16\,\bigl(\Vert\gc^{(2)}\Vert^2\,\dot\gc+\gc^{(3)}\bigr)\,h^2+\go(h^2)
$$
$$  \gN_1(h):=\gv_1(h)-\frac{\gv_1(h)\bullet\gv_0(h)}{\Vert\gv_0(h)\Vert^2}\,\gv_0(h)=2\gc^{(2)}\,h-2\bigl(\Vert\gc^{(2)}\Vert^2\dot\gc+\gc^{(3)}\bigr)\,h^2+\go(h^2)
 $$
$$ \gn_1(h):=\frac{\gN_1(h)}{\Vert\gN_1(h)\Vert}= \\
\frac{\gc^{(2)}}{\Vert\gc^{(2)}\Vert}+\Bigl( -\Vert\gc^{(2)}\Vert\,\dot\gc+\frac{\gc^{(3)}\bullet\gc^{(2)}}{\Vert\gc^{(2)}\Vert^3}\,\gc^{(2)}-\frac 1{\Vert\gc^{(2)}\Vert}\,\gc^{(3)} \Bigr)\,h+\go(h)$$
$$ \ba{rl} \gN_2(h):= & \ds \gv_2(h)-\frac{\gv_2(h)\bullet\gv_0(h)}{\Vert\gv_0(h)\Vert^2}\,\gv_0(h)- \frac{\gv_2(h)\bullet\gn_1(h)}{\Vert\gn_1(h)\Vert^2}\,\gn_1(h) \\
= & \ds 4\,\Bigl( \Vert\gc^{(2)}\Vert^2\dot\gc
-\frac{\gc^{(3)}\bullet\gc^{(2)}}{\Vert\gc^{(2)}\Vert^2}\,\gc^{(2)}
+\gc^{(3)}\Bigr)\,h^2+\go(h^2)=4\gc^{(3)\perp} h^2+\go(h^2)\ea $$
$$ \gn_2(h):=\frac{\gN_2(h)}{\Vert \gN_2(h)\Vert}=\frac{\gc^{(3)\perp}}{\Vert \gc^{(3)\perp} \Vert} +\go(h^0)\,. $$
\ep
\bpf The third order expansions of $\gc$ at $s$ give $\ds \gv_0(h)=\dot\gc+\frac{\gc^{(3)}}6\,h^2+\go(h^2)$ and
$$ \gv_1(h)=-\dot\gc+2\gc^{(2)}\,h-\frac{13}6\,{\gc^{(3)}}\,h^2+\go(h^2)
\,,\quad \gv_2(h)=\dot\gc+2\gc^{(2)}\,h+\frac{13}6\,{\gc^{(3)}}\,h^2+\go(h^2)\,.$$
Whence the formula for $\gt(h)$ follows as
$$\ba{c} \ds \Vert\gv_0(h)\Vert^2=1-\frac{\Vert\gc^{(2)}\Vert^2}3\,h^2+o(h^2)\,,\quad \Vert\gv_0(h)\Vert^{-2}=1+\frac{\Vert\gc^{(2)}\Vert^2}3\,h^2+o(h^2)\,, \\ \ds \Vert\gv_0(h)\Vert^{-1}=1+\frac{\Vert\gc^{(2)}\Vert^2}6\,h^2+o(h^2)\,. \ea $$
We also have
$$ \gv_1(h)\bullet\gv_0(h)=-1+\frac 73\,\Vert\gc^{(2)}\Vert^2h^2+o(h^2) $$
and hence
$$ \frac{\gv_1(h)\bullet\gv_0(h)}{\Vert\gv_0(h)\Vert^2}=-1+2\Vert\gc^{(2)}\Vert^2h^2+o(h^2)$$
that implies the formula for $\gN_1(h)$. We similarly get:
$$ \ba{rl} \ds \Vert\gN_1(h)\Vert^2= & \ds 4\Vert\gc^{(2)}\Vert^2h^2\,\Bigl( 1- 2\frac{\gc^{(3)}\bullet\gc^{(2)}}{\Vert\gc^{(2)}\Vert^2}\,h+o(h)\Bigr)\,, \\
\ds \Vert\gN_1(h)\Vert^{-2}= & \ds \frac 1{4\Vert\gc^{(2)}\Vert^2h^2}\,\Bigl( 1+ 2\frac{\gc^{(3)}\bullet\gc^{(2)}}{\Vert\gc^{(2)}\Vert^2}\,h+o(h)\Bigr)\,, \\
\ds \Vert\gN_1(h)\Vert^{-1}= & \ds \frac 1{2\Vert\gc^{(2)}\Vert h}\,\Bigl( 1 +\frac{\gc^{(3)}\bullet\gc^{(2)}}{\Vert\gc^{(2)}\Vert^2}\,h+o(h)\Bigr) \ea
 $$
that yields the formula for $\gn_1(h)$. Moreover, in order to compute $\gN_2(h)$, we check:
$$ \gv_2(h)\bullet\gv_0(h)=1-\frac 73\,\Vert\gc^{(2)}\Vert^2h^2+o(h^2)\,,\qquad
 \frac{\gv_2(h)\bullet\gv_0(h)}{\Vert\gv_0(h)\Vert^2}=1-2\Vert\gc^{(2)}\Vert^2h^2+o(h^2)$$
and hence
$$ -\frac{\gv_2(h)\bullet\gv_0(h)}{\Vert\gv_0(h)\Vert^2}\,\gv_0(h)= -\dot\gc+\Bigl(2\Vert\gc^{(2)}\Vert^2\,\dot\gc-\frac {1}6\,\gc^{(3)}\Bigr)\,h^2 +\go(h^2)\,. $$
Furthermore,
$$ \gv_2(h)\bullet\gn_1(h)=4\Vert\gc^{(2)}\Vert^2h^2+o(h^2)\,, $$
$$
 \frac{ \gv_2(h)\bullet\gn_1(h) }{\Vert\gn_1(h)\Vert^2}= \frac 1{h^2}\,(h^2+o(h^2))\,\Bigl( 1 +2\,\frac{\gc^{(3)}\bullet\gc^{(2)}}{\Vert\gc^{(2)}\Vert^2}\,h+o(h)\Bigr)$$
that gives
$$ - \frac{ \gv_2(h)\bullet\gn_1(h)}{\Vert\gn_1(h)\Vert^2}\,\gn_1(h) = -2\gc^{(2)}\,h+\Bigl(2\Vert\gc^{(2)}\Vert^2\dot\gc
-4\frac{\gc^{(3)}\bullet\gc^{(2)}}{\Vert\gc^{(2)}\Vert^2}\,\gc^{(2)}
+2\gc^{(3)}\Bigr)\,h^2+\go(h^2)\,. $$
Putting the terms together, we obtain the expression for $\gN_2(h)$, whereas the formula for $\gn_2(h)$ readily follows. \epf
{\large\sc The case of high codimension}. In case of high codimension $N\geq 3$, we wish to extend the previous result to the higher normals.
For this purpose, we introduce the following
\bdf\label{Dst} Let $\gc:[a,b]\to\gR^{N+1}$ be an open rectifiable curve parameterized by arc-length.
Let $j\in\{1,\ldots,N\}$. The curve $\gc$ is said to be {\em smoothly turning at order $j+1$}, if $\gc$ is of class $C^{j+2}$ and at any point $s\in[a,b]$
the vectors
$(\dot\gc(s),\gc^{(2)}(s),\ldots,\gc^{(j+1)}(s))$ are linearly independent. When $j=N$, the curve is said to be {\em smoothly turning}.
\edf
\br If the curve $\gc$ is closed, the same condition is required at any $s\in\gR$, once the curve is extended by periodicity.
\er
%
\par If a curve is smoothly turning, by choosing $s\in]a,b[$, and omitting to write the dependence on $s$, we set:
\beq\label{frame} \ba{c} \ds \gt=\gn_0:=\dot\gc\,,\qquad \gn_1:=\frac{\gc^{(2)}}{\Vert\gc^{(2)}\Vert}\,, \\
\ds \gc^{(j+1)\perp}:= \gc^{(j+1)} - \sum_{k=0}^{j-1}\bigl(\gc^{(j+1)}\bullet\gn_{k}\bigr)\,\gn_k\,,\quad \gn_{j}:=
\frac{\gc^{(j+1)\perp}}{\Vert \gc^{(j+1)\perp} \Vert}\,,\qquad j=2,\ldots,N\,. \ea \eeq
The Jordan frame $(\gt,\gn_1,\ldots,\gn_{N})$ of the curve $\gc$ at the point $\gc(s)$ satisfies the system:
\beq\label{FSN}
 \dot\gt=\gk_1\,\gn_1\,,\quad \dot\gn_1=-\gk_1\,\gt+\gk_2\,\gn_2\,,\quad \dot\gn_j=-\gk_{j}\,\gn_{j-1}+\gk_{j+1}\,\gn_{j+1}\,,\quad j=2,\ldots,N-1  \eeq
where $\gk_j$ is the $j$-{th} curvature of the curve at $\gc(s)$.
\br\label{RFS} The last equation $\dot\gn_{N}=-\gk_{N}\,\gn_{N-1}$ holds true since the curve $\gc$ is dif\-fe\-ren\-ti\-able $(N+2)$-times at the point $s$.
When $N=2$, it reduces to the third Frenet-Serret equation, $\dot\gb=-\gtau\,\gn$.
Since moreover the vectors $(\dot\gc(s),\gc^{(2)}(s),\ldots,\gc^{(N+1)}(s))$ are linearly independent, the last curvature $\gk_N$ is always non-zero.
\er
\br If the curve $\gc$ is smoothly turning at order $j+1$, where $j<N$, only the first $j+1$ Jordan formulas in \eqref{FSN} are satisfied. \er

\par Following the notation from \eqref{vh}, for $k=0,1,\ldots,N$ and for $h>0$ small we define:
\beq\label{vjh} \gv_k(h):= \left\{ \ba{ll} \ds\frac {\gc(s+(k+1)h)-\gc(s+(k-1)h)}{2h} & \text{if $k$ is even} \\
\ds \frac {\gc(s-(k+2)h)-\gc(s-k h)}{2h} & \text{if $k$ is odd}\,. \ea \right. \eeq
%
%
By performing the Gram-Schmidt procedure to $(\gv_0(h),\gv_1(h),\ldots,\gv_{N}(h))$, we also denote as before
$$ \gt(h)=\gn_0(h):=\frac{\gv_0(h)}{\Vert\gv_0(h)\Vert}\,,\quad
 \gN_1(h):=\gv_1(h)-{\bigl(\gv_1(h)\bullet\gt(h)\bigr)}\,\gt(h)\,,\quad
 \gn_1(h):=\frac{\gN_1(h)}{\Vert\gN_1(h)\Vert} $$
and for $j=2,\ldots,N$
$$ \gN_j(h):= \ds \gv_j(h)-\sum_{k=0}^{j-1} {\bigl(\gv_j(h)\bullet\gn_k(h)\bigr)}\,\gn_k(h)
\,,\quad \gn_j(h):=\frac{\gN_j(h)}{\Vert\gN_j(h)\Vert}\,. $$
\par By using a projection argument, we thus obtain:
\bp\label{PGMN} Let $\gc$ be a smoothly turning curve as in Definition~$\ref{Dst}$, and let $(\gt,\gn_1,\ldots,\gn_{N})$ denote the Jordan frame of $\gc$ at a given point $s\in]a,b[$,
see \eqref{frame}.
%
%
Then we have:
$$ \gt(h)=\gt+\go(1)\,,\quad \gn_j(h)=\gn_j+\go(1) \qquad\fa\,j=1,\ldots,N\,. $$
\ep
\bpf
One clearly has $\gt(h)=\gt+\go(1)$. The first step of the Gram-Schmidt procedure, that yields to the formula of $\gn_1(h)$, actually does not depend on the codimension $N\geq 1$, as soon as the higher derivatives $\gc^{(k)}$, for $k\geq 3$, are not involved. Therefore,
since in $\gR^2$ we clearly have $\gn_1(h)=\gn_1+\go(1)$, the same formula holds true in any codimension $N\geq 2$.
\par In a similar way, the second step of the Gram-Schmidt procedure, that yields to the formula of $\gn_2(h)$, does not depend on the codimension $N\geq 2$, as soon as the higher derivatives $\gc^{(k)}$, for $k\geq 4$, are not involved. Therefore,
since in $\gR^3$ we have $\gn_2(h)=\gt(h)\tim\gn_1(h)$, we get $\gn_2(h)=\gn_2+\go(1)$, and hence the same formula holds true in any codimension $N\geq 3$.
\par If $N=3$, we have $\gn_3(h)=\ast(\gt(h)\wdg\gn_1(h)\wdg\gn_2(h))$, where $\ast$ is the Hodge operator in $\gR^4$.
Moreover, $\ast(\gt\wdg\gn_1\wdg\gn_2)=\pm\gn_3$, according to the orientation of the basis $(\gt,\gn_1,\gn_2,\gn_3)$.
By our choice in \eqref{vjh}, this yields that $ \gn_3(h)=\gn_3+\go(1)$, and the projection argument previously described implies that the same formula holds true for $N\geq 4$.
The assertion is proved by proceeding the same way. \epf
\par In general, the higher order coefficients of the expansions of the terms $\gn_j(h)$ actually depend on the choice of the vectors $\gv_k(h)$ we made in \eqref{vjh},
and their existence in general requires more regularity on the curve $\gc$.
For the sake of completeness, in the appendix we provide the following computation in codimension $N=3$, that extends Proposition~\ref{PGM3}.
\bp\label{PGM4} Let $\gc$ be a smoothly turning curve as in Definition~$\ref{Dst}$, where $N=3$.
%
Then at any the given point $s\in]a,b[$ we have:
\beq\label{th}
\gt(h)=\dot \gt+\frac 16\,\bigl(\Vert\gc^{(2)}\Vert^2\,\dot\gc+\gc^{(3)}\bigr)\,h^2+\go(h^3)\,;
\eeq
\beq\label{N1h} \gn_1(h) = \gn_1+
\Bigl( -\Vert\gc^{(2)}\Vert\,\dot\gc+\frac{\gc^{(3)}\bullet\gc^{(2)}}{\Vert\gc^{(2)}\Vert^3}\,\gc^{(2)}-\frac 1{\Vert\gc^{(2)}\Vert}\,\gc^{(3)} \Bigr)\,h+
{\bf d}\,h^2+\go(h^2)  \eeq
for some vector ${\bf d}$ depending on the values of $\dot\gc,\gc^{(2)},\gc^{(3)}$, and $\gc^{(4)}$ at $s$, see \eqref{d} and \eqref{Om};
\beq\label{N2h}
\gn_2(h)=\gn_2+\frac{{\bf D}}{\Vert\gc^{(3)\perp}\Vert}\,h+\go(h) \eeq
where
\beq\label{D} {\bf D}:= \Bigl(\frac{\Vert\gc^{(3)}\Vert^2}{\Vert\gc^{(2)}\Vert}-\Vert\gc^{(2)}\Vert^3 -
\frac{\bigl(\gc^{(3)}\bullet\gc^{(2)}\bigr)^2}{\Vert\gc^{(2)}\Vert^3} \Bigr)\,\gn_1 +
\frac{\gc^{(3)}\bullet\gc^{(2)}}{\Vert\gc^{(2)}\Vert^2\,\Vert\gc^{(3)\perp}\Vert }\, \gn_2  \eeq
and finally
\beq\label{N3h} \gn_3(h)=\gn_3+\go(h) \,.\eeq
%
\ep
\section{Discrete normals to polygonal curves}\label{Sec:discr}
In this section, we introduce a suitable notion of $j$-th normal indicatrix for polygonals.
In fact, the Gram-Schmidt procedure analyzed in the previous section allows us to prove that for smoothly turning curves, one can find a sequence of inscribed polygonals with infinitesimal mesh such that
the length of their $j$-th normal indicatrix converges to the length of the $j$-th normal $\gn_j$ of $\gc$, see
Theorem~\ref{Tappr}.
\par We first fix some notation
and recall some well-known facts, compare e.g. \cite{Su_curv} for further details.
\par Let $P$ denote an oriented polygonal curve in $\gR^{N+1}$, where $N\geq 2$, with ordered (and non-trivial) segments $\{\s_i\mid i=1,\ldots,m\}$, and let $v_i:=\s_i/\calL(\s_i)$ denote the unit vector corresponding to the oriented segment $\s_i$, so that $v_i\in\SN$ for each $i=1,\ldots,m$, where $\SN$ is the Gauss sphere.
The {\em mesh} of the polygonal is defined by $\mesh P:=\sup\{\calL(\s_i)\mid i=1,\ldots,m \}$.
%
\par Following Milnor \cite{Mi}, the tantrix of $P$ is the curve $\gt_P$ in $\SN$ obtained by connecting $v_i$ with $v_{i+1}$ by a minimal geodesic arc, for each $i$, and its length $\calL_{\SN}(\gt_P)$
agrees with the sum of the turning angles, whence with the {\em total curvature} $\TC(P)$ of $P$.
Moreover, if $P$ and $P'$ are polygonal curves in $\gR^{N+1}$, where $P$ is obtained by replacing a segment $\s$ of $P'$ with the two segments joining the end points of $\s$ with a new vertex, then:
$$\calL(P')\leq \calL(P)\,,\qquad \TC(P')\leq\TC(P)\,. $$
\par Similarly to the length, the {\em total curvature} of a curve $\gc$ in $\gR^{N+1}$ is defined by
\beq\label{TC} \TC(\gc):=\sup\{\TC(P)\mid P \prec \gc\} \eeq
where the supremum is taken among all the polygonal curves $P$ inscribed in $\gc$, say $P \prec \gc$.
\par Let $\gc$ be a rectifiable curve with finite total curvature, $\calL(\gc)+\TC(\gc)<\i$.
Due to the previous monotonicity formulas, a continuity argument yields that for any sequence $\{P_n\}$ of inscribed polygonals satisfying $\mesh P_n\to 0$, one has
$\calL(P_n) \to \calL(\gc)$ and $\TC(P_n)\to\TC(\gc)$ as $n\to\i$.
\par In addition, if $\gc:[a,b]\to\gR^{N+1}$ is parameterized by arc-length, so that $\calL(\gc)=b-a$,
then $\gc$ is Lipschitz-continuous, hence it is differentiable a.e., by Rademacher's theorem. Moreover, the tantrix $\gt:=\dot\gc$
is a function of bounded variation in $\BV((a,b),\gR^{N+1})$ taking values in the Gauss sphere $\SN$, and the essential variation $\Var_\SN(\gt)$ of $\gt$ in $\SN$ agrees with the total curvature $\TC(\gc)$.
Therefore, if $\gc$ is of class $C^1$ one has $\TC(\gc)=\int_a^b\Vert \dot\gt(s)\Vert\,ds$, where $\Vert \dot\gt(s)\Vert=\gk_1(s)$, the first curvature of $\gc$.
We refer to \cite{AFP} for the basic notions concerning one-dimensional $\BV$-functions.
\adl\par\noindent
{\large\sc Projective spaces.} The variation of the $j$-th normal to a smooth curve deals with the directions of the osculating spaces of dimension $j$ and $j+1$ through the curvatures $\gk_j$ and $\gk_{j+1}$. Therefore, we compute distances in the projective space $\RPN$, that is defined by the quotient
$\RPN:=\SN/\sim$, the equivalence relation being
$y\sim \wid y\iff y=\wid y$ or $y=-\wid y$, whence the elements of $\RPN$ are denoted by $[y]$.
The projective space $\RPN$ is naturally equipped with the induced metric
$$d_{\RPN}([y],[\wid y]):=\min\{d_{\SN}(y,\wid y),d_{\SN}(y,-\wid y)\}\,. $$
Similarly to $(\SN,d_{\SN})$, the metric space
$(\RPN,d_{\RPN})$ is complete, and the projection map
$\Pi:\SN\to\RPN$ such that $\Pi(y):=[y]$ is continuous.
Moreover, by the {\em
lifting theorem} it turns out that for any continuous function $u:I\to\RPN$ defined on an interval $I\sb\gR$, there
are exactly two continuous functions $v_i:I\to\SN$ such that
$[v_i]:=\Pi\circ v_i=u$, for $i=1,2$, with $v_2(t)=-v_1(t)$ for
every $t\in I$.
\adl\par\noindent
{\large\sc Discrete normals.} Let $P$ be a polygonal curve as above, and assume that $P$ does not lay in a line segment of $\gR^{N+1}$.
For any $i=1,\ldots,m$, we let $v_i^1$ denote the first unit vector $v_h$, with $h>i$, such that $[v_{h}]\neq [v_i]$, so that the linearly independent vectors $(v_i, v_i^1)$ span a 2-dimensional vector space $\Pi^2(P,{v_i})$, that may be called the {\em discrete osculating $2$-space} of $P$ at $v_i$. We then choose the orthogonal direction to $v_i^1$ in $\Pi^2(P,{v_i})$. Therefore, by the Gram-Schmidt procedure, we let
$$ \gN_1(P,i):=v_i-\bigl(v_i\bullet v_i^1\bigr)\,v_i^1\,, \quad \gn_1(P,i):=\frac{\gN_1(P,i)}{\Vert \gN_1(P,i)\Vert} $$
and consider the equivalence class $[\gn_1(P,i)]$.
If $P$ is closed, we trivially extend the notation by listing the vectors $v_i$ in a cyclical way.
If $P$ is not closed and for some $i>1$ there are no vectors $v_h$, with $h>i$, such that $[v_{h}]\neq [v_i]$, we let $[\gn_1(P,i)]:=[\gn_1(P,i-1)]$.
%
%
\smallskip\par In a similar way, if $N\geq 3$, we now define the {\em discrete $j$-th normal} of $P$, for each $j=2,\ldots,N-1$.
We thus assume that $P$ does not lay in an affine subspace of $\gR^{N+1}$ of dimension lower than $j+1$.
For any $i$, we choose $v_i^1$ as above. By iteration on $k=2,\ldots,j$, once we have defined $v_i^{k-1}=v_l$, we let $v_i^k$
denote the first unit vector $v_h$, with $h>l$, such that $v_i^1,v_i^2,\ldots,v_i^k$ are linearly independent.
Therefore, the vectors $(v_i,v_i^1,v_i^2,\ldots,v_i^{j})$ span a $(j+1)$-dimensional vector space $\Pi^{j+1}(P,{v_i})$,
that may be called the {\em discrete osculating $(j+1)$-space} of $P$ at $v_i$.
\par
By means of the Gram-Schmidt procedure, we then choose the {\em orthogonal direction $\gn_j(P,i)\in\SN$ to
$(v_i^1,v_i^2,\ldots,v_i^{j})$ in $\Pi^{j+1}(P,{v_i})$}, and consider the equivalence class $[\gn_j(P,i)]$.
If $P$ is closed, we trivially extend the notation by listing the vectors $v_i$ in a cyclical way.
If $P$ is not closed and for some $i>1$ there are no $j$ vectors satisfying the linear independence as above,
we let $[\gn_j(P,i)]:=[\gn_j(P,i-1)]$.
\smallskip\par Finally, assume now that $P$  does not lay in an affine subspace of $\gR^{N+1}$ of dimension lower than $N$.
The last discrete normal $[\gn_N(P,i)]$ is given by the equivalence class of the {\em orthogonal directions to the
discrete osculating $N$-space $\Pi^{N}(P,{v_i})$} of $P$ at $v_i$.
\bdf\label{Dnj} With the previous notation, for any $j=1,\ldots,N$, we call {\em discrete $j$-th normal} of $P$ the
curve $[\gn_j](P)$ in $\RPN$ obtained by connecting $[\gn_j(P,i)]$ with $[\gn_j(P,i+1)]$ by means of a minimal geodesic arc in $\RPN$, for each $i=1,\ldots,m$, and also $[\gn_j(P,m)]$ with $[\gn_j(P,1)]$, if $P$ is closed.
\edf
\br When $N=2$, i.e., for polygonal curves in $\gR^3$, the last discrete normal $[\gn_2](P)$ agrees with the discrete binormal analyzed in
\cite{Pe,MStor}. As a consequence, its length agrees with the {\em total absolute torsion} $\TAT(P)$ of the polygonal, namely:
\beq\label{TAT} \calL_{\RP}([\gn_2](P))=\TAT(P)\,. \eeq
\par On the other hand, the first discrete normal $[\gn_1](P)$ is different from the weak normal that we introduced \cite{MStor},
where we exploited the polarity in the Gauss sphere $\SP$. \er
{\large\sc A density result.} The following convergence result implies that our notion of
$j$-th normal to a polygonal curve $P$ is the discrete counterpart of the $j$-th normal to a smooth curve $\gc$.
\bt\label{Tappr} Let $\gc:[a,b]\to\gR^{N+1}$, where $N\geq 2$, be a smoothly turning curve at order $j+1$, for some
$j\in\{1,\ldots,N\}$, see Definition~$\ref{Dst}$.
%
Then there exists a sequence $\{P_n\}$ of inscribed polygonals, with $\mesh P_n\to 0$, such that the length $\calL_\RPN([\gn_j](P_n))$ of the discrete $j$-th normal to $P_n$ converges to the length $\calL_\SN(\gn_j)$ of the $j$-th normal $\gn_j$ to the curve $\gc$, i.e.,
$$ \lim_{n\to\i}\calL_\RPN([\gn_j](P_n))=\int_a^b\Vert\dot\gn_j(s)\Vert\,ds\,. $$
%
%
\et
\br We recall that by the Jordan formulas \eqref{FSN}, for each $s\in]a,b[$ one has
$$ \Vert\dot\gn_j(s)\Vert=\sqrt{\gk_{j}^2(s)+\gk_{j+1}^2(s) }$$
if $j<N$, whereas
$ \Vert\dot\gn_N(s)\Vert=|\gk_N(s)|$ when $j=N$, for the last normal.
Moreover, when $N=2$, the last normal $\gn_2$ and curvature $\gk_2$ agree with the binormal and torsion of the curve $\gc$ in $\gR^3$, respectively. \er
\bpff{\em of Theorem~\ref{Tappr}:} If the curve is not closed, we first extend $\gc$ to a smoothly turning curve at order $j+1$ and defined on a closed interval $[\wid a,\wid b]$ such that $\wid a < a<b<\wid b$
\par
For each $n\in\Nat^+$ large, we consider the polygonal curve $P_n$ inscribed in $\gc$ obtained by connecting the consecutive points $\gc(s^n_i)$, where $s^n_i=a+(b-a)\,i/n$, for $i=0,\ldots, n$, whence $\mesh P_n\to 0$ as $n\to \i$, by the uniform continuity of $\gc$. Arguing in a way very similar to the proof of Proposition~\ref{PGMN} and Proposition~\ref{PGM4},
we infer that for each $n$
\beq\label{dnji} [\gn_j(P_n,i)]=[\gn_j(s^n_i)+\ga_j(s^n_i)\,n^{-1}+\go(n^{-1})]\qquad\fa\,i=1,\ldots,n \eeq
where, we recall, $\Vert\go(n^{-1})\Vert=o(n^{-1})$, and $\ga_j(s)$ is a given $\gR^{N+1}$-valued polynomial only depending on the vectors $\dot\gc(s)$, $\gc^{(2)}(s)$,\ldots, $\gc^{(j+1)}(s)$.
\par
Now, since $\gc$ is of class $C^{j+2}$, by the mean value theorem for each $i>1$ we estimate
$$\Vert\ga_j(s^n_{i-1})-\ga_j(s^n_i)\Vert\leq K\cdot n^{-1} $$ for some real constant $K$ depending on the uniform norm of the vector derivatives $\dot\gc(s)$, $\gc^{(2)}(s)$,\ldots, $\gc^{(j+1)}(s)$, whence definitely on $\gc$. Therefore,
for $n$ large enough so that one has $d_\SN(\gn_j(s^n_{i-1}),\gn_j(s^n_{i}))<\p/2$ for each $i$, by the triangular inequality in $\SN$ we can estimate:
$$ \calL_\RPN([\gn_j](P_n))=\sum_{i=2}^n d_\SN(\gn_j(s^n_{i-1}),\gn_j(s^n_{i}))+o(n^{-1}) $$
where $o(n^{-1})\to 0$ as $n\to\i$. Moreover, viewing the points $\{\gn_j(s^n_{i})\mid i=1,\ldots,n\}$ as the vertices of a polygonal $P_n^j$ of $\SN$ inscribed in $\gn_j$,
since $\mesh P_n^j\to 0$, we get $\calL_\SN(P_n^j)\to\calL_\SN(\gn_j)$ as $n\to\i$, whereas
$$ \calL_\SN(P_n^j)=\sum_{i=2}^n d_\SN(\gn_j(s^n_{i-1}),\gn_j(s^n_{i}))\,,\quad  \calL_\SN(\gn_j)=\int_a^b\Vert\dot\gn_j(s)\Vert\,ds\,. $$
The assertion readily follows. 
\epff
\section{Total curvature estimates for the discrete normals}\label{Sec:est}
In this section, we discuss an upper bound for the total curvature of the last normal to a polygonal curve, Proposition~\ref{PineqN}.
In order to extend the upper bound to the intermediate discrete normals, we shall make use of a projection argument and of suitable integral-geometric formulas
for polygonal curves in $\RPN$, that are obtained by extending the integral-geometric formulas for the length and the geodesic rotation of polygonal curves in $\SN$ due to Alexandrov-Reshetnyak \cite{AR}.
\adl\par\noindent
{\large\sc The last normal.} Let $\gc$ be a smoothly turning curve as in Definition~\ref{Dst}, so
that the equation $\dot\gn_{N}=-\gk_{N}\,\gn_{N-1}$ of the Jordan system for the last normal $\gn_N$ holds, where the
last curvature $\gk_N$ is always non-zero.
If $\gT$ denotes the unit tangent vector to the curve $\gn_N$ in $\SN$, one has
$\gT=-\gn_{N-1}$, whence by \eqref{FSN} we get $|\dot\gT|=\sqrt{\gk^2_{N-1}+\gk^2_N}$ and hence the total curvature of $\gn_N$ is equal to the length of the $(N-1)$-th normal:
$$ \TC(\gn_N)=\calL(\gn_{N-1})=\int_a^b\sqrt{\gk^2_{N-1}(s)+\gk^2_N(s)}\,ds\,.
$$
If e.g. $N=2$, then $\gn_2=\gb$, $\gn_1=\gn$, $\gk_1=\gk$, and $\gk_2=\gtau$, and we thus get:
$$ \TC(\gb)=\calL(\gn)=\int_a^b\sqrt{\gk^2(s)+\gtau^2(s)}\,ds\,.
$$
\par We now prove an analogous inequality concerning the discrete last curvature, that goes back to \cite{MStor} for the case of the discrete binormal to polygonal curves is $\gR^3$.
\bp\label{PineqN} Assume $N\geq 2$. Let $P$ be a polygonal curve in $\gR^{N+1}$ that  does not lay in an affine subspace of $\gR^{N+1}$ of dimension lower than $N$, and let $[\gn_j](P)$ denote the discrete $j$-th normal to $P$, see Definition~$\ref{Dnj}$.
Then we have:
$$ \TC([\gn_N](P))\leq \calL_{\RPN}([\gn_{N-1}](P))+\calL_{\RPN}([\gn_{N}](P))\,. $$
\ep
\bpf Recalling the definition of discrete osculating $N$-space $\Pi^{N}(P,{v_i})$ of $P$ at $v_i$, we defined the discrete normal $[\gn_{N-1}(P,i)]$ of $P$ at $v_i$ as the equivalence class in $\RPN$ of the orthogonal directions to
$(v_i^1,v_i^2,\ldots,v_i^{N-1})$ in $\Pi^{N}(P,{v_i})$, and the last discrete normal $[\gn_N(P,i)]$ as the equivalence class of the orthogonal directions to $\Pi^{N}(P,{v_i})$.
\par
If two consecutive osculating $N$-spaces $\Pi^{N}(P,{v_i})$ and $\Pi^{N}(P,{v_{i+1}})$ are different, otherwise $[\gn_N(P,i)]=[\gn_N(P,i+1)]$, and $\g_i$ is the geodesic arc in $\RPN$ connecting the consecutive points $[\gn_N(P,i)]$ and $[\gn_N(P,i+1)]$ of the last discrete normal $[\gn_N](P)$, then $\g_i$ belongs to the great circle corresponding to the 2-dimensional vector space spanned by the independent vectors $\gn_N(P,i)$ and $\gn_N(P,i+1)$.
\par Assuming also without loss of generality that the osculating $N$-spaces $\Pi^{N}(P,{v_{i+1}})$ and $\Pi^{N}(P,{v_{i+2}})$ are different, too, so that the corresponding geodesic arc $\g_{i+1}$ is non-trivial, too,
then the turning angle between $\g_{i}$ and $\g_{i+1}$ is bounded by the length of the geodesic arc in $\RPN$ connecting the consecutive discrete normals
$[\gn_{N-1}(P,i+1)]$ and $[\gn_{N-1}(P,i+2)]$.
\par This property yields that the sum of the turning angles between the consecutive geodesic arcs of
$[\gn_N](P)$ is bounded by the length $\calL_{\RPN}([\gn_{N-1}](P))$ of $[\gn_{N-1}](P)$, whereas the sum of the curvatures of the geodesic arcs $\g_i$ is equal to the length $\calL_{\RPN}([\gn_{N}](P))$ of $[\gn_{N}](P)$, as required.
\epf
\br\label{RTC} If $N=1$, for a polygonal curve $P$ in $\gR^2$ we clearly have 
$$\calL_{ {{\mathbb{RP}}^1} }([\gn_1](P)) \leq \TC(P)\,, \qquad \TC([\gn_1](P))\leq \TC(P)+\calL_{ {{\mathbb{RP}}^1} }([\gn_1](P))\,. $$
\er
\smallskip\par\noindent
{\large\sc Integral-geometric formulas.} For $0\leq j\leq N-1$ integer, denote by $G_{j+1}\gR^{N+1}$ the Grassmannian of
the unoriented ${(j+1)}$-planes in $\gR^{N+1}$. It is a compact group, and it can be
equipped with a unique rotationally invariant probability measure, that will be denoted by
$\m_{j+1}$. For $p\in G_{j+1}\gR^{N+1}$, we denote by $\p_p$ the
orthogonal projection of $\gR^{N+1}$ onto $p$.
\bex\label{ETC} If $\gc$ is a (rectifiable) curve in $\gR^{N+1}$, the following integral-geometric formula for the length holds true for any $j=0,\ldots, N-1$:
$$ \calL(\gc)=\frac{\s_j}{\s_N}\cdot\int_{G_{j+1}\gR^{N+1}}\calL(\p_p(\gc))\,d\m_{j+1}(p) $$
where $\s_j$ and $\s_N$ are positive constants only depending on $j$ and $N$, respectively, see e.g. \cite[Sec. 4.8]{AR}.
\par Let us also recall the average result due to F\'ary
\cite{Fa}, see e.g. \cite[Prop.~4.1]{Su_curv} for a proof, who showed that the total curvature
of a curve (with finite total curvature) is the average of the total curvatures of all its
projections onto $(j+1)$-planes:
%
\beq\label{igTC}
 \TC(\gc)=\int_{G_{j+1}\gR^{N+1}}\TC(\p_p(\gc))\,d\m_{j+1}(p)\qquad \fa j=0,\ldots, N-1\,. \eeq \eex
\par We now deal with polygonal curves in the sphere $\SN$ and in the projective space $\RPN$. Following \cite{AR},
we denote by $\y_p(x)$ the nearest point to $x$ on the $j$-dimensional sphere $\Sjp:=\SN\cap p$.
It is well-defined by
\beq\label{proj} \y_p(x):=\frac{\p_p(x)}{|\p_p(x)|} \eeq
provided that $x\in\SN$ is not orthogonal to the $(j+1)$-plane $p$, i.e., if $x$ does not belong to the $(N-j-1)$-sphere $\Sjp^\perp$ of $\SN$ given by the {\em polar} to $\Sjp$.
Therefore, if $\g$ is a polygonal curve in $\SN$, it turns out that the projected curve $\y_p(\g)$ is well-defined for $\m_{j+1}$-a.e. $p\in G_{j+1}\gR^{N+1}$.
\par The {\em geodesic rotation} $\gK_g(\g)$ of a polygonal curve $\g$ in $\SN$ is given by the sum of the turning angles at the edges of $\g$, see \cite{AR},
so that clearly $\TC(\g)=\calL_\SN(\g)+\gK_g(\g)$.
%
%
The following integral-geometric formulas, that are proved in \cite[Thm. 6.2.2, p. 190]{AR} for $j=1$, actually hold true for larger ranges of values of $j$.
\bt\label{TAR}
Given a polygonal curve $\g$ in $\SN$, for any integer
$1\leq j\leq N-1$ one has
$$
 \ba{rl}\calL_{\SN}(\g)= & \ds \int_{G_{j+1}\gR^{N+1}}\calL_{\Sjp}(\y_p(\g))\,d\m_{j+1}(p) \\
\gK_g(\g)= & \ds \int_{G_{j+1}\gR^{N+1}}\gK_g(\y_p(\g))\,d\m_{j+1}(p)\,. \ea $$
\et
\bpf Assume $j>1$. For $\m_{j+1}$-a.e. $p\in G_{j+1}\gR^{N+1}$, the cited integral-geometric formula from \cite{AR} implies that the length of the projected curve $\calL_{\Sjp}(\y_p(\g))$ is equal to the averaged integral of the projection of the curve $\y_p(\g)$ onto the unit circles corresponding to the 2-planes $q$
of $\gR^{N+1}$ that are contained in $p$, i.e.,
$$\calL_{\Sjp}(\y_p(\g)) = \int_{G_{2}\gR^{j+1}_p} \calL(\y^p_q(\y_p(\g))) \,d\m^p_{2}(q) $$
where $\m^p_{2}$ is the probability measure corresponding to the Grassmannian $G_{2}\gR^{j+1}_p$, with $\gR^{j+1}_p=p$, and $\y^p_q$ is
the nearest point projection from $\Sjp$ onto the 1-circle $\Sjp\cap q$.
Therefore, we have:
$$ \int\lm_{G_{j+1}\gR^{N+1}}\calL_{\Sjp}(\y_p(\g))\,d\m_{j+1}(p) = \int\lm_{G_{j+1}\gR^{N+1}}\Bigl( \int\lm_{G_{2}\gR^{j+1}_p} \calL(\y^p_q(\y_p(\g))) \,d\m^p_{2}(q) \Bigr) \,d\m_{j+1}(p)=:I\,. $$
Moreover, the iterated integral $I$ on the right-hand side is equal to
$$ I=\int_{G_{2}\gR^{N+1}}\calL_{{\mathbb S}^2_r}(\y_r(\g))\,d\m_{2}(r)  $$
and hence, by applying again the formula from \cite{AR}, we get $I=\calL_{\SN}(\g)$, as required. The formula for the
geodesic rotation $\gK_g(\g)$, when $j>1$, is obtained in a similar way from the case $j=1$.
\epf
\par As a consequence, since $\TC(\y_p(\g))=\calL_{\Sjp}(\y_p(\g))+\gK_g(\y_p(\g))$, one also gets:
$$ \TC(\g)=\int_{G_{j+1}\gR^{N+1}}\TC(\y_p(\g))\,d\m_{j+1}(p)\,. $$
\par Now, denote by $\RPjp$ the projective $j$-space corresponding to the $j$-sphere $\Sjp$, for any $p\in G_{j+1}\gR^{N+1}$, and let $\wid\y_p$ denote the nearest point projection of $\RPN$ onto $\RPjp$, i.e., $\wid\y_p([x]):=[\y_p(x)]$, for $x\in \SN\sm\Sjp^\perp$, where $\y_p$ is given by \eqref{proj}.
Following the proof of Theorem~\ref{TAR}, one similarly obtains:
\bp\label{PAR} Given a polygonal curve $\g$ in $\RPN$, for any integer
$1\leq j\leq N-1$ we have
$$ \ba{rl}\calL_{\RPN}(\g)= & \ds \int_{G_{j+1}\gR^{N+1}}\calL_{\RPjp}(\wid\y_p(\g))\,d\m_{j+1}(p) \\
\gK_g(\g)= & \ds \int_{G_{j+1}\gR^{N+1}}\gK_g(\wid\y_p(\g))\,d\m_{j+1}(p) \ea $$
and hence
$$ \TC(\g)=\int_{G_{j+1}\gR^{N+1}}\TC(\wid\y_p(\g))\,d\m_{j+1}(p)\,. $$
\ep
{\large\sc Projection of normals.} We will also make use of the following
\bp\label{Pproj} Let $P$ be a polygonal curve in $\gR^{N+1}$, where $N\geq 2$. For any $j=1,\ldots,N-1$ and for $\m_{j+1}$-a.e. $p\in G_{j+1}\gR^{N+1}$ we have:
$$ [\gn_j](\p_p(P))= \wid\y_p( [\gn_j](P))\,. $$
For $2\leq j\leq N-1$, we also have
$$ [\gn_{j-1}](\p_p(P))= \wid\y_p( [\gn_{j-1}](P))\,. $$
\ep
\bpf Let $\wid\gn_j$ denote the unit vector corresponding by normalization to the projection $\p_p(\gn_j)$ of a vector $\gn_j$ obtained (as in our definition of discrete $j$-th normal from Sec.~\ref{Sec:discr}) by means of the Gram-Schmidt procedure in $\gR^{N+1}$ to a family $v_1,\dots,v_{j+1}$ of independent vectors.
A part the $\m_{j+1}$-negligible case of degeneracy, it turns out that the point $[\wid\gn_j]\in \RPjp$ agrees with the equivalence class of the unit vector obtained by applying the analogous Gram-Schmidt procedure in $p\in G_{j+1}\gR^{N+1}$ to the projected vectors
$\p_p(v_1),\dots,\p_p(v_{j+1})$. Therefore, the first formula readily follows on account of Definition~\ref{Dnj}, and the second one is proved in a similar way. \epf
\par By Propositions~\ref{PAR} and~\ref{Pproj}, we readily obtain the following
\bc\label{CigP} If $P$ is a polygonal curve in $\gR^{N+1}$, for any $1\leq j\leq N-1$ we have:
%
%
$$ \calL_{\RPN}([\gn_j](P))= \ds \int_{G_{j+1}\gR^{N+1}}\calL_{\RPjp}( [\gn_j](\p_p(P)) )\,d\m_{j+1}(p)\,. $$
\ec
\par In the case $j=1$, we also infer:
\bp\label{PigP} If $P$ is a polygonal curve in $\gR^{N+1}$, where $N\geq 2$, we have:
$$ \calL_{\RPN}([\gn_1](P))\leq\TC(P)\,.
 $$
\ep
\bpf By Remark~\ref{RTC}, for $\m_2$-a.e. $p\in G_2\gR^{N+1}$ one has $\calL_{\RPjp}( [\gn_j](\p_p(P)) )\leq \TC(\p_p(P))$. Therefore, the inequality follows from Corollary~\ref{CigP} and from the integral-geometric formula \eqref{igTC} for the total curvature, by monotonicity of the averaged integral.
\epf
%
%
{\large\sc The intermediate normals.} Finally, by using Propositions~\ref{PAR} and~\ref{Pproj}, we are able to extend the total curvature estimate to the intermediate normals.
\bp\label{Pineqj} Let $P$ be a polygonal curve in $\gR^{N+1}$, where $N\geq 2$, and let $[\gn_j](P)$ denote the discrete $j$-th normal to $P$, see Definition~$\ref{Dnj}$.
Then for every $j=2,\ldots,N$ we have:
$$ \TC([\gn_j](P))\leq \calL_{\RPN}([\gn_{j-1}](P))+\calL_{\RPN}([\gn_{j}](P))\,. $$
Moreover, for $j=1$ we have
$$\TC([\gn_1](P)) \leq \calL_{\SN}(\gt_P)+\calL_{\RPN}([\gn_{1}](P))\,,\qquad \calL_{\SN}(\gt_P)=\TC(P)\,. $$
\ep
\bpf If $j=N$, the assertion follows from Proposition~\ref{PineqN}. If $N\geq 3$ and $j=2,\ldots,N-1$, by Proposition~\ref{PAR} we have
$$ \TC([\gn_j](P))=\int_{G_{j+1}\gR^{N+1}}\TC(\wid\y_p([\gn_j](P)))\,d\m_{j+1}(p)\,. $$
Therefore, by Proposition~\ref{Pproj} we can write:
$$ \TC([\gn_j](P))=\int_{G_{j+1}\gR^{N+1}}\TC([\gn_j](\p_p(P)))\,d\m_{j+1}(p)\,. $$
\par By applying Proposition~\ref{PineqN}, with $j$ instead of $N$, to the last curvature of $\p_p(P)$, we have
$$ \TC([\gn_j](\p_p(P)))\leq \calL_{\RPjp}([\gn_{j-1}](\p_p(P)))+\calL_{\RPjp}([\gn_{j}](\p_p(P))) $$
for $\m_{j+1}$-a.e. $p\in G_{j+1}\gR^{N+1}$, so that again by Proposition~\ref{Pproj} we get:
$$ \TC([\gn_j](\p_p(P)))\leq \calL_{\RPjp}( \wid\y_p( [\gn_{j-1}](P)) )+\calL_{\RPjp}( \wid\y_p( [\gn_{j}](P))) $$
and hence, by the monotonicity of the averaged integral,
$$ \TC([\gn_j](P))\leq \int_{G_{j+1}\gR^{N+1}}\bigl[
\calL_{\RPjp}( \wid\y_p( [\gn_{j-1}](P)) )+\calL_{\RPjp}( \wid\y_p( [\gn_{j}](P))) \bigr]  \,d\m_{j+1}(p)\,.$$
\par
By applying again the integral-geometric formulas from Proposition~\ref{PAR}, we get:
$$ \ba{r} \ds \int_{G_{j+1}\gR^{N+1}}\bigl[\calL_{\RPjp}([\gn_{j-1}](\p_p(P)))+\calL_{\RPjp}([\gn_{j}](\p_p(P)))\bigr]\,d\m_{j+1}(p) \\ \ds = \calL_{\RPN}([\gn_{j-1}](P))+\calL_{\RPN}([\gn_{j}](P)) \ea $$
and the claim readily follows.
Finally, the case $j=1$ follows from Remark~\ref{RTC}, by means of a similar argument.
\epf
\section{The relaxed total variation of the normals to a curve}\label{Sec:rel}
In this section, we introduce a relaxed notion of total variation of the $j$-th normal to a curve.
Due to the lack of monotonicity, we are led to follow the approach introduced by Alexandrov-Reshetnyak \cite{AR}, that involves the notion of {\em modulus}.
\br\label{Rmon} Differently to what happens for length and total curvature, the monotonicity formula fails to hold in general for the length of the discrete
$j$-th normal to polygonals.
More precisely, if $P$ and $P'$ are polygonal curves in $\gR^{N+1}$, where $P$ is obtained by replacing a segment $\s$ of $P'$ with the two segments joining the
end points of $\s$ with a new vertex, then it may happen that $\calL_{\RPN}([\gn_j](P'))>\calL_{\RPN}([\gn_j](P))$ for some $j=1,\ldots,N$.
This feature was observed in \cite{MStor} concerning the length of the discrete binormal to polygonal curves in $\gR^3$, i.e., about the functional $P\mapsto\calL_{\RP}([\gn_2](P))$, that agrees with our notion of total absolute torsion $\TAT(P)$ of the polygonal, see \eqref{TAT}.
\er
\bex\label{Emon} (cf. \cite{MStor}). Let $P$ be a polygonal made of six segments $\s_i$, for $i=1,\ldots,6$, where the first three ones and the last three ones lay on
two different planes $\Pi_1$ and $\Pi_2$.
Then the tantrix $\gt_P$ connects with geodesic arcs in $\SP$ the consecutive points $v_i:=\s_i/\calL(\s_i)$, for $i=1,\ldots, 6$, where the triplets
$v_1,v_2,v_3$ and $v_4,v_5,v_6$ lay on two geodesic arcs, which are inscribed in the great circles corresponding to the vector spaces spanning the planes $\Pi_1$ and $\Pi_2$, respectively.
If both the angles $\a$ and $\be$ of $\gt_P$ at the points $v_3$ and $v_4$ are small, then $\TAT(P)=\a+\be$.
\par
\begin{figure}
	\centering
	\includegraphics[width=0.50\textwidth]{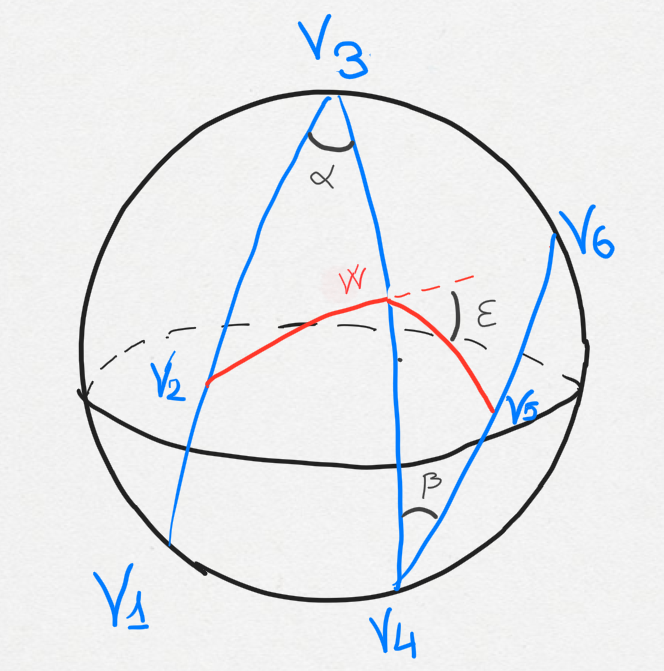}
		\caption{The tantrix of the polygonal $P$, in blue color, and of the inscribed polygonal $P'$, in red color.
\newline
The drawing is courtesy offered by the young artist Sofia Saracco.}
\label{fig:controex}
\end{figure}

 Let $P'$ be the inscribed polygonal obtained by replacing the segments $\s_3$ and $\s_4$ of $P$ with the segment $\s$ between the first point of $\s_3$ and the last point of $\s_4$.
The tantrix $\gt_{P'}$ connects with geodesic arcs the consecutive points $v_1,v_2,w,v_5,v_6$, where the point $w$ lays in the minimal geodesic arc between $v_3$ and $v_4$.
Now, assume that the turning angle $\e$ of $\gt_{P'}$ at the point $v_5$ satisfies $\a<\e<\p/2$, and that the two geodesic triangles with vertices $v_2,v_3,w$ and $w,v_4,v_5$ have the same area. By suitably choosing the position of the involved vertices, and by using the Gauss-Bonnet theorem in the computation, it turns out that $\TAT(P')-\TAT(P)=2(\e-\a)>0$, see Figure~\ref{fig:controex}.
\eex
%
%
\par We recall that the {\em modulus} $\m_\gc(P)$ of a polygonal curve $P$ inscribed in a curve $\gc$ of $\gR^{N+1}$ is the maximum of the diameter
of the arcs of $\gc$ determined by the consecutive vertices in $P$.
\par We correspondingly notice that, if $\gc$ is a polygonal curve itself, there exists $\e>0$ such that any polygonal $P$ inscribed in $\gc$ and with modulus $\m_c(P)<\e$ satisfies
$\gt_P=\gt_\gc$, whence $[\gn_j](P)=[\gn_j](\gc)$ for each $j=1,\ldots,N$.
It suffices indeed to take $\e$ lower than half of the mesh of the polygonal $\gc$, so that in every segment of $\gc$ there are at least two vertices of $P$.
\smallskip\par The above facts motivate us to introduce the following:
\bdf\label{DFj} Let $\gc$ be a curve in $\gR^{N+1}$. The {\em relaxed total variation of the $j$-th normal} to $\gc$ is given by
\beq\label{Fjc} \F_j(\gc):=\lim_{\e\to 0^+}\sup\{\calL_\RPN([\gn_j](P)) \mid P \prec \gc\,,\,\,\m_c(P)<\e\} \qquad j=1,\ldots,N \eeq
where $[\gn_j](P)$ is the discrete $j$-th normal to the inscribed polygonal
$P$, see Definition~\ref{Dnj}.
\edf
\par By the previous remark, in fact, for any polygonal curve $P$ in $\gR^{N+1}$ we have
\beq\label{P2}
 \F_j(P) = \calL_\RPN([\gn_j](P)) \qquad\fa\,j=1,\ldots,N\,. \eeq
We can thus re-write the integral-geometric formulas for polygonals in Corollary~\ref{CigP} as:
\beq\label{igP}
\F_j(P) = \ds \int_{G_{j+1}\gR^{N+1}}\F_j(\p_p(P))\,d\m_{j+1}(p)\,,\qquad 1\leq j\leq N-1\,. \eeq
\br\label{Rigp} For future use, we point out that when $j>1$ one similarly gets
$$\F_{j-1}(P) = \ds \int_{G_{j+1}\gR^{N+1}}\F_{j-1}(\p_p(P))\,d\m_{j+1}(p)\,. $$ \er
\br
%
%
When $N=2$, according to \eqref{TAT}, it turns out that the relaxed total variation of the last normal agrees with the notion of total absolute torsion for curves $\gc$ in $\gR^3$ that we analyzed in \cite{MStor}, namely
$$\F_2(\gc)=\TAT(\gc)\,. $$
\par Notice that, in order to extend formula \eqref{igP} to the relaxed total variation of the normals to a curve $\gc$,
we cannot argue as for the total curvature, see Example~\ref{ETC}, where one applies the
monotone convergence theorem to a sequence of approximating polygonals with $P_n\prec P_{n+1}\prec\gc$ for each $n$, compare e.g. \cite[Prop.~4.1]{Su_curv}.
In fact, we have seen in Remark~\ref{Rmon} that the monotonicity property fails to hold in this context.
\er
{\large\sc Properties.}
If $\F_j(\gc)<\i$ for some $j=1,\ldots,N$, for any sequence $\{P_n\}$ of polygonal curves inscribed in $\gc$ and satisfying $\m_\gc(P_n)\to 0$, one has $\sup_n\calL_\RPN([\gn_j](P_n))<\i$.
Also, one can find an optimal sequence as above in such a way that $\calL_\RPN([\gn_j](P_n))\to\F_j(\gc)$ as $n\to\i$.
\smallskip\par Moreover, the relaxed total variation of the first normal is always lower than the total curvature:
\bp\label{PF1TC} For any curve $\gc$ in $\gR^{N+1}$, according to formula \eqref{TC}, we have
\beq\label{F1TC} \F_1(\gc)\leq\TC(\gc)\,.
\eeq \ep
\bpf If $\TC(\gc)<\i$, the following result from \cite[Thm.~2.1.3]{AR} holds true: for each $\e>0$ there exists $\delta>0$ such that if $\gamma$ is an arc of $\gc$ with geodesic
diameter lower than $\delta$, the length of $\g$ is smaller than $\e$.
As a consequence, if $\gc$ has finite total curvature, one has:
$$\TC(\gc)=\lim_{\e\to 0^+}\sup\{\TC(P) \mid P \prec \gc\,,\,\,\m_c(P)<\e\}\,. $$
Therefore, inequality \eqref{F1TC} readily follows from Proposition~\ref{PigP}. \epf
\br In general, the strict inequality holds in \eqref{F1TC}. In fact, for e.g. a polygonal curve $P$ in $\gR^2$, in the quantity
%
$\calL_{ {{\mathbb{RP}}^1}}([\gn_1](P))$ we take distances in the projective line, so that a contribution of $\TC(P)$ given by a turning angle $\t$ greater than $\p/2$, corresponds to a contribution $\p-\t$ for the length of $[\gn_1](P)$.
\er
\par As a consequence of Theorem~\ref{Tjcurv}, we readily obtain the following continuity property.
\bp\label{PFjc} Let $N\geq 2$ and $\gc$ be a curve in $\gR^{N+1}$ such that $\F_{j}(\gc)<\i$ and $\F_{j-1}(\gc)<\i$ for some $j=2,\ldots,N$.
Then, for any sequence $\{P_n\}$ of inscribed polygonals satisfying $\m_\gc(P_n)\to 0$ one has:
$$ \lim_{n\to\i} \calL_\RPN([\gn_j](P_n))=\F_j(\gc)\,. $$
In the case $j=1$, the same conclusion holds true for any curve $\gc$ satisfying $\TC(\gc)<\i$.
\ep
\par Therefore, for smoothly turning curves, the following explicit formulas for the relaxed total variation of the normals hold:
\bp\label{PFjcsmooth}
Let $\gc:[a,b]\to\gR^{N+1}$, where $N\geq 2$, be a smoothly turning curve at order $j+1$, for some
$j\in\{1,\ldots,N\}$, see Definition~$\ref{Dst}$.
Then we have
$$ \F_j(\gc)=\int_a^b\Vert\dot\gn_j(s)\Vert\,ds $$
where, we recall, $\Vert\dot\gn_j(s)\Vert=\sqrt{\gk_j^2(s)+\gk_{j+1}^2(s)}$, when $j<N$, and
$\Vert\dot\gn_N(s)\Vert=|\gk_N(s)|$, when $j=N$.
%
%
\ep
\bpf By the density theorem~\ref{Tappr}, the hypotheses of Theorem~\ref{Tjcurv} are clearly satisfied.
Therefore, the assertions follow from Proposition~\ref{PFjc}, on account of the Jordan formulas \eqref{FSN}, and of Remark~\ref{RFS} in the case $j=N$. \epf
\section{Weak normals to a non-smooth curve}\label{Sec:weakn}
In this section, we analyze a weak notion of $j$-th normal to a curve $\gc$ in $\gR^{N+1}$ such that $\F_{j}(\gc)<\i$.
We are able to define a Lipschitz-continuous curve $[\gn_j](\gc)$ on $\RPN$, parameterized by arc-length and satisfying
\beq\label{Varnj}\calL_{\RPN}([\gn_j](\gc))=\F_{j}(\gc) \eeq
in such a way that for any sequence of inscribed polygonals converging to $\gc$, the length of the discrete $j$-th normals converges to the length of the curve $[\gn_j](\gc)$.
%
%
\par We shall make use of arguments taken from \cite{MStor} for the case of the binormal indicatrix of curves in $\gR^3$. Since the compactness argument
relies on the curvature estimates for polygonals from Proposition~\ref{Pineqj}, we need to assume in addition that $\F_{j-1}(\gc)<\i$, in the case $j>1$, and that the curve $\gc$ has finite total curvature, when $j=1$.
\bt\label{Tjcurv} Let $N\geq 2$ and $\gc$ be a curve in $\gR^{N+1}$ such that $\F_{j}(\gc)<\i$ and $\F_{j-1}(\gc)<\i$ for some $j=2,\ldots,N$.
There exists a rectifiable curve $[\gn_j](\gc):[0,L_j]\to\RPN$ parameterized by arc-length, where $L_j:=\F_{j}(\gc)$, so that \eqref{Varnj} holds true, satisfying the following property.
For any sequence $\{P_n\}$ of inscribed polygonal curves, let $\bfg^j_n:[0,L_j]\to\RPN$ denote for each $n$ the parameterization with constant velocity
of the discrete $j$-th normal  $[\gn_j](P_n)$ to $P_n$, see Definition~$\ref{Dnj}$.
If $\m_\gc(P_n)\to 0$, then $\bfg^j_n\to [\gn_j](\gc)$ uniformly on $[0,L_j]$ and
$$\calL_{\RPN}(\bfg^j_n)=\calL_{\RPN}([\gn_j](P_n))\to\calL_{\RPN}([\gn_j](\gc))$$
as $n\to\i$, where, we recall, $\calL_{\RPN}([\gn_j](P_n))=\F_j(P_n)$.
Moreover, the arc-length derivative of the curve $[\gn_j](\gc)$ is a function of bounded variation.
Finally, in the case $j=1$, for any curve
$\gc$ in $\gR^{N+1}$ satisfying $\TC(\gc)<\i$, one has $\F_1(\gc)<\i$ and the same conclusion as above holds true.
\et
%
%
\par
It is quite easy to construct a smooth  curve whose $(j-1)$-th curvature is infinite, while its $j$-th curvature is finite, or even zero: it is enough to take the curve in an affine space of the appropriate dimension. We show an explicit example of a rectifiable curve in $\gR^3$ whose curvature is infinite, while its torsion is zero.
\bex Let $\g:[0,1]\to\gR^3$ be defined as follows
$$\g(t)=\left\{
\begin{array}{ll}(0,0,0) & \ \text{if}\  t=0\\ (t^2\sin(2\pi/t),t^2\cos(2\pi/t),0) & \ \text{if}\   t\neq0\,.
\end{array}\right.$$
When $t$ ranges from $1/(n+1)$ to $1/n$ the curve makes a complete loop around the origin at a distance lower than $1/n^2$, for each $n\in\Nat^+$. Therefore, the curve is of finite length, since its length may be estimated with the convergent sum $\sum_n 2\pi/n^2$, while its total curvature is infinite. Its torsion is obviously zero, since the curve is planar.

One may add a small non-planarity to the example, e.g. making the last coordinate be $e^{-1/t^2}$ instead of zero, causing the torsion to be bigger than zero, but still finite, and still having the curvature infinite.
\eex
\par Motivated by Theorem~\ref{Tjcurv}, that will be proved below, we introduce the following
\bdf\label{Dnjc} Under the hypotheses of Theorem~$\ref{Tjcurv}$, the curve $[\gn_j](\gc)$ is called {\em weak $j$-th normal} to the curve $\gc$.
\edf
\par We also notice that Proposition~\ref{PFjc} is a direct
consequence of Theorem~\ref{Tjcurv}. Finally, at the end of this section we also prove the validity of the following integral-geometric formula:
\bc\label{Cigc} For curves $\gc$ in $\gR^{N+1}$ satisfying $\F_{j}(\gc)+\F_{j-1}(\gc)<\i$ for some integer $2\leq j <N$, we have:
\beq\label{igc}\F_j(\gc)=\int_{G_{j+1}\gR^{N+1}}\F_j(\p_p(\gc)))\,d\m_{j+1}(p)\,. \eeq
When $j=1$, the same formula holds true for curves $\gc$ in $\gR^{N+1}$ satisfying $\TC(\gc)<\i$.
\ec
\adl\par\noindent
\bpff{\em of Theorem~\ref{Tjcurv}:}
It is divided into eight steps. When $j>1$, in Steps~1-2, we obtain the curve $[\gn_j](\gc)$ by means of an optimal approximating sequence. In Steps~3-4, where we exploit the polarity of the last normal, we deal with the case $j=N$.
In Steps~5-7, where we first make use of the integral-geometric formula \eqref{igP} for polygonals, we analyze the case $1<j<N$ of the intermediate normals.
Finally, in Step~8 we deal with the case $j=1$ of the first normal.
\smallskip\par\noindent{\em Step~1:} Assume $j>1$. Choose an optimal sequence $\{P_n\}$ of polygonal curves inscribed in $\gc$ such that $\m_\gc(P_n)\to 0$ and
$L^j_n\to L_j$, where $L^j_n:=\calL_{\RPN}([\gn_j](P_n))$, the curve $[\gn_j](P_n)$ being the discrete $j$-th normal to $P_n$, see Definition~\ref{Dnj},
and, we recall, $L_j:=\F_j(\gc)$.
If $\F_j(\gc)=0$, the proof is trivial.
Assuming $0<\F_j(\gc)<\i$, for $n$ large enough so that $L^j_n>0$, we also denote by $[\gn_j](P_n):[0,L^j_n]\to\RPN$
the arc-length parameterization of the curve $[\gn_j](P_n)$.
%
%
%
\par Define $\bfg^j_n:[0,L_j]\to\RPN$ by $\bfg^j_n(s):=[\gn_j](P_n)((L^j_n/L_j)s)$,
so that $\Vert \dot\bfg^j_n(s)\Vert = L^j_n/L_j$ a.e., where $L^j_n/L_j\to 1$. By Ascoli-Arzela's
theorem, we can find a (not relabeled) subsequence of $\{\bfg^j_n\}$ that uniformly
converges in $[0,L_j]$ to some Lipschitz continuous function $\bfg^j:[0,L_j]\to\RPN$. Whence, $\bfg^j$ is differentiable a.e., by Rademacher's theorem, whereas
by lower-semicontinuity $\Vert\dot\bfg^j(s)\Vert\leq 1$ for a.e. $s\in[0,L_j]$.
\smallskip\par\noindent{\em Step~2:} We claim that $\dot \bfg^j_n\to \dot \bfg^j$ strongly in $L^1$.
As a consequence, we deduce that $\Vert \dot \bfg^j\Vert =1$ a.e., and hence, denoting
$\bfg^j=\gn_j[\gc]$, that
$$ \calL_{\RPN}(\gn_j[\gc])=\int_0^{L_j}\Vert \dot\bfg^j(s)\Vert \,ds=L_j=\F_{j}(\gc)\,. $$
\par In order to prove the claim, in this step we choose a (not relabeled) continuous lifting of the curve $\bfg^j$, so that $\bfg^j:[0,L_j]\to\SN$, and for $n$ large enough, we identify the curve $\bfg^j_n$ with its (not relabeled) continuous lifting $\bfg^j_n:[0,L_j]\to\SN$ such that $d_\SN(\bfg^j_n(0),\bfg^j(0))<\p/2$.
Consider the tantrix $\tt^j_n(s):=\dot\bfg^j_n(s)/\Vert \dot\bfg^j_n(s)\Vert$ of the curve $\bfg^j_n$, where, we recall,
$\Vert \dot\bfg^j_n(s)\Vert = L^j_n/L_j$ a.e., with $L^j_n/L_j\to 1$. We have
$\Var(\tt^j_n)\leq\TC(\bfg^j_n)$, whereas by Proposition~\ref{Pineqj}, we can estimate the total curvature of each curve $\bfg^j_n$ as follows:
$$ \TC(\bfg^j_n)=\TC([\gn_j](P_n))\leq \calL_{\RPN}([\gn_{j-1}](P_n))+\calL_{\RPN}([\gn_{j}](P_n))\,. $$
Since we assumed $\F_{j-1}(\gc)<\i$, we also have $\sup_n\calL_{\RPN}([\gn_{j-1}](P_n))<\i$, whence we get:
$$ \sup_n\Var(\tt^j_n)\leq\sup_n \TC(\bfg^j_n)<\i\,. $$
As a consequence, by compactness, a further subsequence of $\{\dot\bfg^j_n\}$ converges weakly-* in the $\BV$-sense to some $\BV$-function $v:[0,L_j]\to\gR^{N+1}$.
The claim follows if we show that $v(s)=\dot\bfg^j(s)$ for a.e. $s\in[0,L_j]$. In fact, this property yields that the sequence $\{\dot \bfg^j_n\}$ converges strongly in $L^1$
to the function $\dot \bfg^j$. In particular, by lower semicontinuity it turns out that $\dot\bfg^j$ is a function of bounded variation.
\par Now, using that by
Lipschitz-continuity
$$ \bfg^j_n(s)=\bfg^j_n(0)+\int_0^s \dot\bfg^j_n(\l)\,d\l\qquad \fa\, s\in [0,L_j]
$$
and setting
$$ V(s):=\bfg^j(0)+\int_0^s v(\l)\,d\l\,,\qquad s\in [0,L_j] $$
by the weak-* $\BV$ convergence $\dot\bfg^j_n\wc v$, which implies the
strong $L^1$ convergence, we have $\bfg^j_n\to V$ in $L^\i$,
hence $\dot\bfg^j_n\to \dot V=v$ a.e. on $[0,L_j]$. But we already
know that $\bfg^j_n\to \bfg^j$ in $L^\i$, thus we get $v=\dot\bfg^j$.
\smallskip\par\noindent{\em Step~3:} Assume now $j=N$. Let $\{\wid P_n\}$ denote any sequence of polygonal curves inscribed in $\gc$ such that $\m_\gc(\wid P_n)\to 0$. We show that
possibly passing to a subsequence, the discrete $N$-th normals $[\gn_{N}](\wid P_n)$ uniformly converges (up to reparameterizations, as above)
to the curve $[\gn_{N}](\gc)$.
\par For this purpose, we recall from Sec.~\ref{Sec:discr} that the discrete osculating $N$-space $\Pi^{N}(P,{v_i})$ of a polygonal $P$ at the unit vector $v_i$
is given by the hyperplane spanned by consecutive points in the Gauss sphere $\SN$ which correspond to
consecutive vertexes of the tantrix $\gt_P$.
Moreover, the last discrete normal $[\gn_N(P,i)]$ is identified by the orthogonal directions to $\Pi^{N}(P,{v_i})$,
whence by the {\em polar} in the projective space $\RPN$ to the hyper-sphere corresponding to the discrete osculating $N$-space of $P$ at $v_i$.
\par Now, if $\{P_n\}$ is the optimal sequence of the previous steps (with $j=N$), conditions $\m_\gc(\wid P_n)\to 0$ and $\m_\gc(P_n)\to 0$
yield that the {\em Frech\'et distance} (see e.g. \cite[Sec.~1]{Su_curv}) between the two sequences $\{\gt_{P_n}\}$ and $\{\gt_{\wid P_n}\}$ goes to zero.
Recalling our Definition~\ref{Dnj} of discrete $N$-th normal $[\gn_N](P)$, by the continuity of the Gram-Schmidt procedure and of the polarity transformation, it turns out that
the Frech\'et distance between $[\gn_{N}](\wid P_n)$ and $[\gn_{N}](P_n)$ goes to zero, but we already know that a sub-sequence of $\{[\gn_{N}](P_n)\}$ uniformly converges
to the curve $[\gn_{N}](\gc)$, as required.
\smallskip\par\noindent{\em Step~4:}
If $j=N$ and $\{\wid P_n\}$ is the (not relabeled) subsequence obtained in Step~3, by repeating the argument in Step~1 we infer that
the limit function $\bfg^N$ is unique. As a consequence, a contradiction argument yields that the whole sequence $\{\wid\bfg^N_n\}$ uniformly converges to $\bfg^N$
and that the limit curve $\bfg^N=[\gn_{N}](\gc)$ does not depend on the choice of the sequence $\{\wid P_n\}$ of
inscribed polygonals satisfying $\m_\gc(\wid P_n)\to 0$. Therefore, the curve $[\gn_{N}](\gc)$ is
identified by $\gc$. Arguing as in Step~2, we finally infer that $\calL_{\RPN}([\gn_{N}](\wid P_n))\to\calL_{\RPN}([\gn_{N}](\gc))$.
\smallskip\par\noindent{\em Step~5:} Assume now $1<j<N$. We claim that the function $g_j(p):=\F_j(\p_p(\gc))$, for $p\in G_{j+1}\gR^{N+1}$, belongs to the summable class $L^1(G_{j+1}\gR^{N+1},\m_{j+1})$.
\par In fact, if $\{P_n\}$ is an optimal sequence of inscribed polygonals from Steps~1-2, so that $\F_j(P_n)\to\F_j(\gc)$,
using the integral-geometric formula \eqref{igP}, by Fatou's Lemma we have
$$ \F_j(\gc)\geq \int_{G_{j+1}\gR^{N+1}}\liminf_{n\to\i}\F_j(\p_p(P_n))\,d\m_{j+1}(p)\,. $$
The sequence $\{\p_p(P_n)\}$ of polygonals is inscribed in $\p_p(\gc)$ and satisfies $\m_{\p_p(\gc)}(\p(P_n))\to 0$. Moreover, by the previous inequality, and using Definition~\ref{DFj},
we infer that $\F_j(\p_p(\gc))<\i$ for $\m_{j+1}$-a.e. $p\in G_{j+1}\gR^{N+1}$.
On account of Remark~\ref{Rigp}, we similarly obtain $\F_{j-1}(\p_p(\gc))<\i$ for $\m_{j+1}$-a.e. $p$.
\par
Therefore, by Steps~3-4, where we take $j=N$ (and work with the last discrete normal to the projected curve), we infer that
$\F_j(\p_p(P_n))\to\F_j(\p_p(\gc))=g_j(p)$ for $\m_{j+1}$-a.e. $p$, whence $g_j$ is measurable and
\beq\label{ineq}
\int_{G_{j+1}\gR^{N+1}}g_j(p)\,d\m_{j+1}(p)=\int_{G_{j+1}\gR^{N+1}}\F_j(\p_p(\gc))\,d\m_{j+1}(p)\leq\F_j(\gc)<\i \eeq
so that the claim readily follows.
\smallskip\par\noindent{\em Step~6:} Let $\{\wid P_n\}$ denote any sequence of polygonal curves inscribed in $\gc$ such that $\m_\gc(\wid P_h)\to 0$. We show that
$\F_j(\wid P_n)=\calL_{\RPN}([\gn_{j}](\wid P_n))\to\F_j(\gc)$.
\par In fact, if $\{P_n\}$ is the optimal sequence from the previous step, by \eqref{igP} for each $n$ we estimate
\beq\label{int} |\F_j(\wid P_n)-\F_j(P_n)|\leq \int_{G_{j+1}\gR^{N+1}}|\F_j(\p_p(\wid P_n))-\F_j(\p_p(P_n))|\,d\m_{j+1}(p)\,. \eeq
Moreover, again by Definition~\ref{DFj}, for $\m_{j+1}$-a.e. $p$ we can find $\e(p)>0$ such that if $P\prec \gc$ satisfies $\m_{\p_p(\gc)}(\p_p(P))<\e(p)$, then $\F_j(\p_p(P))<2\,\F_j(\p_p(\gc))$. Also, by compactness of the Grassmannian $G_{j+1}\gR^{N+1}$ we get $\ol\e:=\inf_p\e(p)>0$. Therefore, since
$\m_{\p_p(\gc)}(\p_p(P))\leq\m_{\gc}(P)$, we can find $\ol n$ such that for any $n>\ol n$
$$ |\F_j(\p_p(\wid P_n))-\F_j(\p_p(P_n))|\leq 4\,\F_j(\p_p(\gc))=4\,g_j(p)  $$
for $\m_{j+1}$-a.e. $p\in G_{j+1}\gR^{N+1}$.
Arguing as above, by Step~4, where we take $j=N$, we infer that
$\F_j(\p_p(\wid P_n))\to\F_j(\p_p(\gc))$ and hence that $|\F_j(\p_p(\wid P_n))-\F_j(\p_p(P_n))|\to 0$ for $\m_{j+1}$-a.e. $p$. Since
$g_j\in L^1(G_{j+1}\gR^{N+1},\m_{j+1})$, by dominated convergence the integral in equation \eqref{int} goes to zero as $n\to\i$, whence
$\F_j(\wid P_n)\to\F_j(\gc)$.
\smallskip\par\noindent{\em Step~7:} Now, if $1<j<N$, for any sequence $\{\wid P_n\}$ of inscribed polygonal curves with $\m_\gc(\wid P_n)\to 0$, as in Steps~1-2 we infer that possibly passing to a subsequence
$\wid\bfg^j_n\to \wid\bfg^j$ uniformly on $[0,L_j]$ to some curve $\wid\bfg^j$ parameterized in arc-length.
If $\{P_n\}$ is the optimal sequence, we denote by $\wih P_n$ the polygonal given by the common refinement of $P_n$ and $\wid P_n$. The uniform limit of (a subsequence of) the corresponding sequence $\{\wih\bfg^j_n\}$ is equal to the uniform limit of both $\{\wid\bfg^j_n\}$ and $\{\bfg^j_n\}$.
This yields that $\wid\bfg^j=[\gn_j](\gc)$. Finally, the proof is completed by arguing as in Step~4.
\smallskip\par\noindent{\em Step~8:} In the case $j=1$, the first statement follows from Proposition~\ref{PF1TC}.
The proof proceeds as in the case $j>1$ above, on account of the following straightforward modifications.
Firstly, in Step~2, by Proposition~\ref{Pineqj} we can estimate the total curvature of each curve $\bfg^1_n$ as follows:
$$ \TC(\bfg^1_n)=\TC([\gn_1](P_n))\leq\calL_{\SN}(\gt_{P_n})+\calL_{\RPN}([\gn_{1}](P_n))\,,\qquad \calL_{\SN}(\gt_{P_n})=\TC(P_n) $$
and hence the role of the functional $\F_{j-1}(\cdot)$ is played by the total curvature $\TC(\cdot)$, when $j=1$. In fact, since we assumed $\TC(\gc)<\i$, we also have $\sup_n\TC(P_n)<\i$, whence we get $ \sup_n \TC(\bfg^1_n)<\i$.
Secondly, in Step~5, by using this time the integral-geometric formula \eqref{igTC}, with $j=1$, we infer that $\F_1(\p_p(\gc))<\i$ and $\TC(\p_p(\gc))<\i$ for $\m_{2}$-a.e. $p\in G_{2}\gR^{N+1}$. We omit any further detail.
\epff
\br\label{Rcomp} In Step~2, we could have proved the $L^1$-convergence of $\dot\g^j_n$ to $\dot\g^j$ by applying the
Kolmogorov-Riesz-Frech\'et compactness theorem, thus showing that
$$ \lim_{|\delta|\to 0}\sup_n\int_0^{L_j } |\dot \g^j_n(s+\delta)-\dot \g^j_n(s)|\,ds=0\,. $$
However, for each $s\in[0,L_j]$ and for $\delta>0$ small we can estimate
$$
\int_0^{L_j } |\dot \g^j_n(s+\delta)-\dot \g^j_n(s)|\,ds  \leq c\cdot\delta\cdot \TC(\g^j_n)  $$ 
for some absolute constant $c$ and hence we need the additional assumption  $\F_{j-1}(\gc)<\i$. On the other hand, we showed that $\dot\g^j$ is a function of bounded variation, a property that will be used in Sec. \ref{Sec:curvmeas}, where we introduce the curvature measures by means of the first variation formula of the length of the curve $\bfg^j$, see \eqref{var}.
\er
\bpff{\em of Corollary~\ref{Cigc}:} Since the integral-geometric formula holds true for polygonals, it suffices to argue in a way very similar to Step~6, on account of the dominated convergence theorem.
\epff

\section{Relationship with the smooth normals}\label{Sec:relation}
%
In this section, we wish to find a wider class of smooth curves $\gc$ for which our weak $j$-th normal $[\gn_j](\gc)$ is strictly related to the classical $j$-th normal $\gn_j$ to $\gc$, see Definition~\ref{Dms}. In fact, for smoothly turning curves, see Definition~\ref{Dst}, this property is outlined in
Proposition~\ref{Csmooth}.
As we shall see below, the main property we need to preserve is the {\em existence and continuity of the osculating $(j+1)$-spaces}.
%
%
\adl\par\noindent{\large\sc Smoothly turning curves.} 
As a first consequence of Proposition~\ref{PFjc}, by the density theorem~\ref{Tappr} and the Jordan formulas \eqref{FSN}, in Proposition~\ref{PFjcsmooth} we obtained that the relaxed total variation of the $j$-th normal agrees with the length of the smooth $j$-th normal $\gn_j$.
We now see that the weak $j$-th normal $[\gn_j](\gc)$ is equivalent to the smooth $j$-th normal.
\bp\label{Csmooth} Let $\gc:[a,b]\to\gR^{N+1}$, where $N\geq 2$, be a smoothly turning curve at order $j+1$, for some
$j\in\{1,\ldots,N\}$, see Definition~$\ref{Dst}$.
Then, the weak $j$-th normal $[\gn_j](\gc)$ agrees (up to a lifting from $\RPN$ to $\SN$) with the arc-length parameterization of the smooth $j$-th normal $\gn_j$ to $\gc$. More precisely, if $\Pi:\SN\to\RPN$ is the canonical projection, one has
$$ [\gn_j](\gc)(t)  = \Pi(\gn_j(\psi_j(t)))\qquad \fa\,t\in[0,L_j] $$
where $\psi_j:[0,L_j]\to[a,b]$ is the inverse of the bijective and $C^1$-class
transition function
\beq\label{vfj} \vf_j(s):=\int_a^s\Vert\dot\gn_j(\l)\Vert\,d\l\,,\qquad s\in[a,b] \eeq
and, we recall,
$ L_j:=\F_j(\gc)=\calL_{\RPN}([\gn_j](\gc))$.
\ep
\bpf Going back to the proof of Theorem~\ref{Tappr}, it turns out that
the sequence $\{P_n\}$ of inscribed polygonals satisfies $\m_\gc(P_n)\to 0$. Moreover, formula \eqref{dnji},
where, we recall, the coefficients $\ga_j(s^n_i)$ are equibounded in terms of the
uniform norm in $[a,b]$ of the vector derivatives $\gc^{(k)}$, for $k=1,\ldots,j+1$, implies that the Frech\'et distance between
the curves $[\gn_j](P_n)$ and $\gn_j$ goes to zero as $n\to\i$. Therefore, one has
\beq\label{direct} \Pi(\gn_j(s))=[\gn_j](\gc)(\vf_j(s))\qquad \fa\,s\in[a,b]\,. \eeq
Moreover, the linear independence of the vectors
$\dot\gc(s),\gc^{(2)}(s),\ldots,\gc^{(j+1)}(s)$ for any $s\in[a,b]$, on account of
the Jordan equations \eqref{FSN} and of formulas \eqref{frame}, yields that the arc-length derivative $\dot\gn_j(s)$ is non-zero for every $s$.
The assertion readily follows.
\epf
{\large\sc Milder conditions.}
In our paper \cite{MStor} on curves in $\gR^3$, we noticed that the existence of the osculating plane to a smooth curve $\gc$, is guaranteed by the requirement that at each
point $s$ there exists a non-zero higher order derivative $\gc^{(k)}(s)$.
In fact, by computing the derivatives in the identity $\dot\gc\bullet\dot\gc=0$ one sees that the osculating plane at $\gc(s)$, say $\Pi^2(\gc,s)$, is given by
$\gc(s)+\span\{\dot\gc(s),\gc^{(k)}(s)\}$, where $k$ is the smallest integer $k>1$ such that $\gc^{(k)}(s)\neq 0_{\gR^3}$.
Therefore, the 2-vector $\dot\gc(s)\wdg\gc^{(k)}(s)$ provides an orientation to the osculating plane,
and the unit normal $\gn(s)$ is given by applying the Gram-Schmidt procedure to the couple of vectors $\dot\gc(s),\gc^{(k)}(s)$.
Moreover, it turns out that the second derivative $\gc^{(2)}$ is zero only at a finite set of point,
but in general the normal $\gn(s)$ fails to be continuous when these ones are inflection points.
However, the osculating plane $\Pi^2(\gc,s)$ is a continuous function of the arc-length parameter.
This property ensures that the normal vector $\gn$
(and hence the binormal vector $\gb=\gt\tim\gn$, too) is continuous when seen as a function in the projective plane $\RP$.
The following example of {\em mildly smoothly turning curve}, see Definition~\ref{Dms}, is taken from \cite{MStor}.
\bex\label{Eflex} Let $\gc:[-1,1]\to\gR^3$ be the curve satisfying $\gc(0)=0_{\gR^3}$ and with derivative
$$\dot \gc(s)={1\over\sqrt 2}\,\bigl(1,s^2,\sqrt{1-s^4}\bigr),\quad s\in[-1,1] $$
so that $\Vert\dot \gc(s)\Vert\equiv 1$. We compute
$$\gc^{(2)}(s)={\sqrt 2 s\over \sqrt{1-s^4}}\,\bigl(0,\sqrt{1-s^4},-s^2\bigr)\,,\quad \gc^{(3)}(s)=\sqrt 2\,\Bigl(0,1,{s^2(s^4-3)\over (1-s^4)^{3/2}} \Bigr)\,. $$
Therefore, if $0<|s|<1$ we have $\gc^{(2)}(s)\neq 0_{\gR^3}$ and hence
$$ \gn(s)={s\over |s|}\,\bigl(0,\sqrt{1-s^4},-s^2\bigr)\,,\quad \gb(s)={s\over |s|}\,{1\over\sqrt 2}\,\bigl(-1,s^2,\sqrt{1-s^4}\bigr)\,. $$
Furthermore, for $0<|s|<1$ we get:
$$ \gk(s):=\Vert \gc^{(2)}(s)\Vert={\sqrt 2|s|\over \sqrt{1-s^4}}\,,
\quad\gtau(s):={\bigl(\dot \gc(s)\tim\gc^{(2)}(s)\bigr)\bullet \gc^{(3)}(s) \over \Vert \gc^{(2)}(s)\Vert^2}=-{\sqrt 2s\over \sqrt{1-s^4}} $$
and hence $\gk(s)\to 0$ and $\gtau(s)\to 0$ as $s\to 0$, both $\gk$ and $\gtau$ are summable functions in $L^1(-1,1)$, and the Frenet-Serret formulas hold true separately in the open intervals $]-1,0[$ and $]0,1[$.
\par Since $\gt(0)=2^{-1/2}(1,0,1)$, $\gc^{(2)}(0)=0_{\gR^3}$, and $\gc^{(3)}(0)=2^{-1/2}(0,1,0)$,
the osculating plane at $\gc(0)$ is
$$\Pi^2(\gc,0)=0_{\gR^3}+\span\{2^{-1/2}(1,0,1),\,2^{-1/2}(0,1,0)\}$$
and by the Gram-Schmidt procedure we get
$\gn(0)=(0,1,0)$ and $ \gb(0)=2^{-1/2}\,( -1,0,1)$. Therefore, even if the unit normal and binormal are not continuous at $s=0$,
since $[\gn(s)]\to[\gn(0)]$ and $[\gb(s)]\to[\gb(0)]$ as $s\to 0$, they are both continuous as functions with values in $\RP$. For future use, we finally compute
\beq\label{exdot}  {\dot\gn(s) \over \Vert\dot\gn(s)\Vert}={s\over |s|}\,\bigl(0,\,-s^2,-\sqrt{1-s^4}\bigr)\,, \qquad s\neq 0\,. \eeq
\eex
\par For curves in $\gR^{N+1}$, where $N>2$, the above argument concerning the osculating 2-plane continues to hold.
In order to deal with the high dimension osculating spaces, the analogous sufficient condition is given by the
existence of $j+1$ independent derivatives $\gc^{(k)}(s)$ of the curve near each point $\gc(s)$.
\bdf\label{Dms} Let $\gc:[a,b]\to\gR^{N+1}$, where $N\geq 2$, be an open rectifiable curve parameterized in arc-length.
The curve is said to be {\em mildly smoothly turning at order $j+1$}, where  $j\in\{1,\ldots,N\}$, if for each $s\in[a,b]$ the function $\gc$ is of class $C^m$ in a neighborhood of $s$, for some
integer $m\geq j+2$, and there exist $j$ integers $1< i_2<\ldots < i_{j+1}< m$ such that the $(j+1)$-vector
$(\dot\gc\wdg\gc^{(i_2)}\wdg\cdots\wdg\gc^{(i_{j+1})})(s)$ is non-trivial. When $j=N$, the curve is said to be {\em mildly smoothly turning}.
\edf
\br If the curve $\gc$ is closed, the same condition is required at any $s\in\gR$, once the curve is extended by periodicity.
\er
\par With these assumptions,
%
in fact, the osculating $(j+1)$-space $\Pi^{j+1}(\gc,s)$ to the curve at $\gc(s)$, is spanned by the $(j+1)$-vector obtained by choosing the smallest indexes $i_k$ as above, see formula \eqref{osc},
and it moves continuously along the curve, Proposition~\ref{Posc}.
Moreover, the first $j$ unit normals are defined
by following the idea due to Jordan.
\bdf\label{Dnjmild} Let $\gc$ be a mildly smoothly turning curve at order $j+1$, where $j<N$, and let $1<i_2<\ldots < i_{j+1}$ be the smallest integers such that the $(j+1)$-vector
$(\dot\gc\wdg\gc^{(i_2)}\wdg\cdots\wdg\gc^{(i_{j+1})})(s)$ is non-trivial. The $j$-th normal $\gn_j(s)$ is defined by the last term in the Gram-Schmidt procedure to the ordered list of independent vectors $\dot\gc(s),\,\gc^{(i_2)}(s),\ldots,\gc^{(i_{j+1})}(s)$. If $\gc$ is a mildly smoothly turning curve, we also set
$\gn_N:=\ast(\gt\wdg\gn_1\wdg\cdots\wdg\gn_{N-1})$, where $\ast$ is the Hodge operator in $\gR^{N+1}$.
\edf
\par Of course, a smoothly turning curve at order $j+1$ is mildly smoothly turning at the same order, and the above property at a higher order implies the same one at lower orders.
\par We now show that the features we obtained in the smoothly turning case, can be extended by considering equivalence classes of antipodal points
in the Gauss sphere $\SN$.
\par More precisely, we recover the convergence result, Proposition~\ref{Pappr},
the representation formula for the relaxed functional $\F_j(\gc)$, Proposition~\ref{PFjcmild},
and the relationship between the weak $j$-th normal $[\gn_j](\gc)$ from Theorem~\ref{Tjcurv} and the smooth $j$-th normal, Proposition~\ref{Cmild}.
\par We first notice that if a smooth curve fails to satisfy the linear independence property in Definition~\ref{Dms},
then the osculating $(j+1)$-space fails to be continuous, in general.
%
\bex\label{Eflat} Let $f:\mathbb R\to\mathbb R$ be the $C^\infty$ but not analytic function given by
$$f(x)\ :=\ \left\{\begin{array}{ll}
e^{-1/x^2} & \mbox{if \ $x\neq 0$}\\
 0 & \mbox{if \  $x=0$}\,.
\end{array}\right.$$
The function $f$ has all derivatives vanishing in zero. Let us consider the curve $\g:[-1,1]\to\mathbb R^3$ defined as
$$
\g(t)\ :=\ \left\{ \ba{ll} \bigl(t,\,f(t),\,0\bigr) & \mbox{if \ $t\leq 0$}\\
\bigl(t,\,0,\,f(t)\bigr)& \mbox{if \ $t\geq 0$}\,.
\ea
\right.
$$
The curve $\g$ is smooth ($C^\infty$), but since all its derivatives $\g^{(2)},\dots,\g^{(n)},\dots$ vanish in zero, it does not satisfy the assumptions in Definition~\ref{Dms}. The same is true if one considers a re-parametrization $\gc$ of $\g$ in arc-length.

Since for $t\leq 0$ the curve lies in the plane $\pi_1=\{z=0\}$ and for $t\geq0$ it lies in the plane $\pi_2=\{y=0\}$, the torsion of the curve is always zero, $\gb$ is constant out of $t=0$, and $\gb$ and $\gn$ jump of an angle of $\pi/2$ at $t=0$. By modifying the plane $\pi_2$, it is immediate to find an example in which the curve has both the normal $\gn$ and binormal $\gb$ jumping of an arbitrary angle $\alpha$ at $t=0$. Notice that since $\gt$ is continuous and $\gb=\gt\times\gn$, the jump angle $\alpha$ must be the same for both $\gn$ and $\gb$.

Moreover, the example is easily adapted to curves in spaces of higher dimension having an arbitrary number of normals jumping of arbitrary angles. Notice, though, that since the last normal $\gn_{N}$ is determined by the vectors $\gt,\gn_1,\ldots,\gn_{N-1}$, the angle of jump of the last normal $\gn_{N}$ is determined by those of the other normals.
\eex
{\large\sc Properties.} In the sequel, without loss of generality we deal with open curves, and $j\in\{1,\ldots,N\}$, with $N\geq 3$, if not differently specified.
\bp\label{Pcomp} If $\gc$ is a mildly smoothly turning curve at order $j+1$, there exists a finite set $\SS$ of points in $]a,b[$ such that the $(j+1)$-vector
$(\dot\gc\wdg\gc^{(2)}\wdg\cdots\wdg\gc^{(j+1)})(s)$ is non-trivial on $]a,b[\sm\SS$.
Moreover, the first $j$ formulas in the Jordan system \eqref{FSN} are satisfied in each connected component of $]a,b[\sm \SS$,
and the corresponding curvature terms $\gk_h$ are continuous functions on $]a,b[$, that may possibly be equal to zero only at the
singular points $s_i\in\SS$.
Moreover, if the curve is mildly smoothly turning, the last formula in the Jordan system \eqref{FSN} holds true, too, on $]a,b[\sm\SS$.
\ep
\bpf Since linear independence is an open property, a compactness argument yields the first assertion. The other ones readily follow. \epf
\par
The main feature is the existence and continuity of the osculating $(j+1)$-spaces along the curve.
In fact, equipping the set of {\em unoriented $(j+1)$-planes} with the canonical metric, we have:
\bp\label{Posc} If a curve $\gc$ is mildly smoothly turning at order $j+1$, the osculating $(j+1)$-space $\Pi^{j+1}(\gc,s)$ is well-defined and continuous, as
$s\in]a,b[$.
\ep
\bpf For fixed $s\in]a,b[$, consider the $j+1$ vectors $\gv_k(h)$ given by \eqref{vjh}, for $k=0,\ldots,j$, and let $1<i_2<\ldots < i_{j+1}$ be the smallest integers such that the $(j+1)$-vector
$(\dot\gc\wdg\gc^{(i_2)}\wdg\cdots\wdg\gc^{(i_{j+1})})(s)$ is non-trivial. For $h\neq 0$ small, by using as before the Taylor expansions of $\gc$ centered at $s$ and at order $m$,
and writing the wedge product, we obtain
$$ \gv_0(h)\wdg \gv_1(h)\wdg\cdots \wdg\gv_{j}(h) = \l\,(\dot\gc\wdg\gc^{(i_2)}\wdg\cdots\wdg\gc^{(i_{j+1})})(s)\cdot h^p+\go(h^m) $$
where the integer $p:=(i_2+\ldots+i_{j+1})-j\in\Nat^+$ and the factor $\l$ is a non-zero real number that depends on the indexes $i_k$, through the Taylor expansions. This yields that
$$
\frac { \gv_0(h)\wdg \gv_1(h)\wdg\cdots \wdg\gv_{j}(h) }{| \gv_0(h)\wdg \gv_1(h)\wdg\cdots \wdg\gv_{j}(h)|} =
\Bigl(\frac h{|h|}\Bigl)^p\gu_{j+1}(s)+o(1) $$ 
where
$$ \gu_{j+1}(s):=
\frac{ (\dot\gc\wdg\gc^{(i_2)}\wdg\cdots\wdg\gc^{(i_{j+1})})(s)}{|(\dot\gc\wdg\gc^{(i_2)}\wdg\cdots\wdg\gc^{(i_{j+1})})(s)|}\,. $$
By smoothness, letting $h\to 0$ we infer that the non-zero unit $(j+1)$-vector $\gu_{j+1}(s)$ provides an orientation to the osculating
$(j+1)$-space $\Pi^{j+1}(\gc,s)$ to the curve at $\gc(s)$, and actually
\beq\label{osc} \Pi^{j+1}(\gc,s)=\gc(s)+\span\{ \dot\gc(s),\,\gc^{(i_2)}(s),\ldots,\gc^{(i_{j+1})}(s)\}\,. \eeq
\par Now, by Proposition~\ref{Pcomp} it turns out that for each $s\in]a,b[\sm\SS$
$$ \gu_{j+1}(s)=
\frac{ (\dot\gc\wdg\gc^{(2)}\wdg\cdots\wdg\gc^{({j+1})})(s)}{|(\dot\gc\wdg\gc^{(2)}\wdg\cdots\wdg\gc^{({j+1})})(s)|} $$
where $\SS$ is a finite set, and hence the $(j+1)$-vector function $s\mapsto \gu_{j+1}(s)$ may fail to be continuous at the points $s_i\in\SS$.
%
However, since the smooth vectors $\gv_k(h)$ are defined in terms of Taylor expansions of $\gc$ at $s$, and $\gc$ is of class $C^m$ near each $s_i$, where $m>i_{j+1}$, it turns out that at any point $s_i\in\SS$ one has
$$\gu_{j+1}(s_i)=\pm\gu_{j+1}(s_i-)=\pm \gu_{j+1}(s_i+)\,. $$
\par
Since the topology induced by the canonical metric of unoriented $(j+1)$-spaces is equivalent to the one induced by the equivalence classes of unoriented unit $(j+1)$-vectors, the continuity property follows.
%
%
\epf
\br\label{Rmild}
For smoothly turning curves in the sense of Definition~\ref{Dst}, we always have $i_k=k$ for each $k=2,\ldots,j+1$,
and
the $(j+1)$-vector function $s\mapsto\gu_{j+1}(s)$ is continuous in $]a,b[$, actually of class $C^1$.
\par More generally, if the curve $\gc$ is mildly smoothly turning at order $j+1$,
at each point $s_i\in \SS$ the normals may be discontinuous.
However, denoting by $f(s\pm)$ the right and left limits of a function $f$ at the point $s$, the continuity of the osculating $(j+1)$-space along the curve implies the equalities
$$\gn_k(s_i-)=\pm\gn_k(s_i+) \qquad\fa\,k=1,\ldots,j $$
and hence {\em the first $j$ unit normals are continuous when seen as a
function into the projective space $\RPN$}.
\par Moreover, by our assumptions the $(j+1)$-vector $\gu_{j+1}(s)$ is of class $C^1$ in each connected component of $]a,b[\sm\SS$.
More precisely, it turns out that the osculating $(j+1)$-space function $s\mapsto\Pi^{j+1}(\gc,s)$ is of class $C^1(]a,b[)$, w.r.t. the canonical metric of unoriented $(j+1)$-spaces in $\gR^{N+1}$.
%
In addition, the curvature terms $\gk_{j-1}$ and $\gk_j$ are always non-zero on $]a,b[\sm\SS$. We thus obtain:
\beq\label{jumps}
 {\dot\gn_j(s_i-)\over \Vert\dot\gn_j(s_i-)\Vert}=\pm {\dot\gn_j(s_i+)\over \Vert\dot\gn_j(s_i+)\Vert}\in\SN  \eeq
according to formula \eqref{exdot} from Example~\ref{Eflex}.
\er
\par We now readily extend the convergence result obtained in Theorem~\ref{Tappr}.
\bp\label{Pappr} Let $\gc$ be a mildly smoothly turning curve at order $j+1$, for some $1\leq j\leq N$.
Then there exists a sequence $\{P_n\}$ of inscribed polygonals, with $\mesh P_n\to 0$, such that
$$ \lim_{n\to\i}\calL_\RPN([\gn_j](P_n))=\int_a^b\Vert\dot\gn_j(s)\Vert\,ds\,, $$
where $\gn_j$ is given by Definition~$\ref{Dnjmild}$.
\ep
\bpf If the curve is not closed, we first extend $\gc$ to a mildly smoothly turning curve at order $j+1$ and defined on a closed interval $[\wid a,\wid b]$ such that $\wid a < a<b<\wid b$.
The proof then proceeds in a very similar way to the one of Theorem~\ref{Tappr}. Notice, in fact, that with our assumptions the equalities
\eqref{dnji} continue to hold for each $n$. We omit any further detail. \epf
\par Moreover, the representation formula for the relaxed total variation of the $j$-th normal, see Proposition~\ref{PFjcsmooth}, continues to hold:
\bp\label{PFjcmild} If $\gc$ is a mildly smoothly turning curve at order $j+1$, for some $1\leq j\leq N$, we have
$$ \F_j(\gc)=\int_a^b\Vert\dot\gn_j(s)\Vert\,ds<\i\,. $$
%
\ep
\bpf By Proposition~\ref{Pappr}, the curve $\gc$ satisfies the hypotheses of Theorem~\ref{Tjcurv}.
Therefore, the claim follows from Proposition~\ref{PFjc} and from the Jordan formulas in Proposition~\ref{Pcomp}. \epf
\par Finally, we recover the relationship in Proposition~\ref{Csmooth} between the weak $j$-th normal $[\gn_j](\gc)$ from Theorem~\ref{Tjcurv} and the smooth $j$-th normal.
\bp\label{Cmild} Let $\gc$ be a mildly smoothly turning curve at order $j+1$, and let $\gn_j$ be given by Definition~$\ref{Dnjmild}$. Then we have:
$$ [\gn_j](\gc)(t) = \Pi(\gn_j(\psi_j(t)))\qquad \fa\,t\in[0,L_j] $$
where $\Pi:\SN\to\RPN$ is the canonical projection, $\psi_j:[0,L_j]\to[a,b]$ is the inverse of the bijective and absolutely continuous transition function \eqref{vfj},
and, we recall,
$ L_j:=\F_j(\gc)=\calL_{\RPN}([\gn_j](\gc))$.
\ep
\bpf We argue in a way very similar to the proof of Corollary~\ref{Csmooth}. In fact, in the proof of Proposition~\ref{Tappr},
the sequence $\{P_n\}$ of inscribed polygonals satisfies $\m_\gc(P_n)\to 0$, whereas formula \eqref{dnji}
implies again that the Frech\'et distance between
the curves $[\gn_j](P_n)$ and $\gn_j$ goes to zero as $n\to\i$, so that \eqref{direct} holds true.
This time, by Proposition~\ref{Posc} we deduce that
the arc-length derivative $\dot\gn_j(s)$ of the smooth $j$-th normal in Definition~\ref{Dnjmild} is non-zero for every $s\in]a,b[$ except to a finite set of singular points $s_i$.
This property implies that the transition function \eqref{vfj} is bijective and absolutely continuous, as required. \epf
\section{Curvature measures}\label{Sec:curvmeas}
The {\em curvature force} was introduced in \cite{CFKSW}, see also \cite{Su_curv}, as the distributional derivative of the tangent indicatrix of curves in $\gR^{N+1}$ with finite total curvature, the starting point being the computation of the first variation of the length of the curve.
Using similar arguments, when $N=2$, the {\em torsion force} was discussed in \cite{MStor}, where we considered tangential variations of the length of the tantrix. We now see that similar arguments can be repeated for the weak $j$-th normals.
As before, in the sequel we deal with open curves.
\par To this purpose, we recall that in Theorem~\ref{Tjcurv}, we showed that the arc-length derivative of the curve $[\gn_j](\gc)$ in $\RPN$ is a function of bounded variation. For simplicity, we denote here by $\bfg^j:[0,L_j]\to\SN$ a continuous lifting of the curve $[\gn_j](\gc)$, so that $\dot\bfg^j$ is a function of bounded variation, with $\Vert\dot\bfg^j\Vert\equiv 1$.
Moreover, we have:
$$ \calL_\SN(\bfg^j)=\calL_\RPN([\gn_j](\gc))=\F_j(\gc)\,. $$
\par We assume that $\bfg^j_\e$ is a variation of $\bfg^j$ under which the motion of each point $\bfg^j(t)$ is smooth in time and with initial velocity $\x(t)$, where $\x:[0,L_j]\to\gR^{N+1}$ is a Lipschitz continuous function with $\x(0)=\x(L_j)=0$, so that
$\dot\x(t)$ is defined for a.e. $t$, by Rademacher's theorem.
\par
Denoting by $D\dot\bfg^j$ the finite measure given by the distributional derivative of $\dot\bfg^j$, the first variation formula of the length of the curve $\bfg^j$ gives:
\beq\label{var} \d_\x\calL_{\SN}(\bfg^j):={d\over d\e}\,\calL_\SN(\bfg^j_\e)_{\vert \e=0}=\int_0^{L_j} \dot\bfg^j(t)\bullet\dot\x(t)\,dt=:-\lan D\dot\bfg^j,\x\ran\,. \eeq
{\large\sc The polygonal case}. If $\gc$ is a polygonal curve $P$, the weak $j$-th normal agrees with the discrete $j$-th normal $[\gn_j](P)$ from Definition~\ref{Dnj}, obtained by connecting the consecutive points $[\gn_j(P,i)]$ with minimal geodesic arcs in $\RPN$. Therefore, the arc-length derivative of the lifting $\bfg^j$ has a discontinuity in correspondence eventually to the points $[\gn_j(P,i)]$, where the norm of the jump is equal to the turning angle between
the consecutive geodesic arcs meeting at $[\gn_j(P,i)]$. Therefore, the total variation of the measure $D\dot\g^j$ is equal to the total curvature of the curve
$\dot\g^j$ in $\gR^{N+1}$, and hence to the sum $\calL_{\RPN}([\gn_j](P))+ \TC_{\RPN}([\gn_j](P))$, where $\TC_{\RPN}$ is the intrinsic total curvature of the
curve in $\RPN$. We omit any further detail.
\adl\par\noindent
{\large\sc Smoothly turning curves}. Assume now that the curve $\gc$ is smoothly turning at order $j+1$, Definition~\ref{Dst}. By Proposition~\ref{Csmooth}, possibly considering the antipodal continuous lifted function of $[\gn_j](\gc)$, for every $t\in[0,L_j]$ we have
$ \bfg^j(t)=\gn_j(\psi_j(t))$.
Then, by changing variable $t=\vf_j(s)$ we can write
\beq\label{Dgj} \lan D\dot\bfg^j,\x\ran=-\int_a^{b} \dot\bfg^j(\vf_j(s))\bullet{d\over ds}\,[\x(\vf_j(s))]\,ds \eeq
and hence, using that
\beq\label{dergj} \dot\bfg^j(t)={\dot\gn_j(s)\over \Vert\dot\gn_j(s)\Vert}\,,\qquad t=\vf_j(s)  \eeq
and integrating by parts, since $\x(\vf_j(a))=\x(\vf_j(b))=0$ we obtain:
\beq\label{repr} \lan D\dot\bfg^j,\x\ran=-\int_a^{b} {\dot\gn_j(s)\over \Vert\dot\gn_j(s)\Vert}\bullet{d\over ds}\,[\x(\vf_j(s))]\,ds =
\int_a^{b} {d\over ds}\,{\dot\gn_j(s)\over \Vert\dot\gn_j(s)\Vert}\bullet \x(\vf_j(s))\,ds\,.
\eeq
Therefore, the function $\dot\bfg^j$ is of class $C^1(]a,b[)$, and denoting by $\calL^1$ the Lebesgue measure in $\gR$, it turns out that the distributional derivative of $\dot\g^j$ is an absolutely continuous measure
\beq\label{meas} D\dot\bfg^j=\vf_{j\,\#}\m_j\,,\qquad \m_j:={d\over ds}\,{\dot\gn_j(s)\over \Vert\dot\gn_j(s)\Vert}\,\calL^1\pri]a,b[ \eeq
given by the push forward of the measure $\m_j$ by the function $t=\vf_j(s)$.
%
%
\par In general, when $j<N$ the denominator $\Vert\dot\gn_j\Vert$ in the formula \eqref{dergj} involves two curvatures. Therefore, the explicit computation of the density of the measure $\m_j$ involves five normals and four curvatures.
We now consider in particular the simpler case of the last normal.
\bex\label{Emeas} When \,$j=N$, we recall the last two Jordan formulas:
$$ \dot\gn_{N-1}=-\gk_{N-1}\,\gn_{N-2}+\gtau\,\gn_{N}\,, \qquad \dot\gn_N=-\gtau\,\gn_{N-1} $$
where we have denoted $\gtau:=\gk_N$, the last curvature
(that is, the torsion, when $N=2$, in which case the Frenet-Serret formulas give $\gn_0=\gt$, $\gn_1=\gn$, $\gk_1=\gk$, and $\gn_2=\gb$).
Denoting by $\sgn\gtau$ the constant sign of the non-zero smooth function $\gtau(s)$, we thus obtain:
$$ {\dot\gn_N(s)\over \Vert\dot\gn_N(s)\Vert} =  -\sgn\gtau\cdot \gn_{N-1}(s)\,,\qquad
{d\over ds}\,{\dot\gn_N(s)\over \Vert\dot\gn_N(s)\Vert} =\sgn\gtau\cdot\Bigl( \gk_{N-1}\,\gn_{N-2}-\gtau\,\gn_{N}\Bigr)(s)\,.
$$
%
%
\par Now, we restrict to consider tangential variations in formula \eqref{var}, i.e., {\em we assume in addition that $\x(t)\in T_{\bfg^j(t)}\SN$ for each $t$}.
We correspondingly deduce that the {\em tangential component} $D^\top\bfg^N$ of the measure $D\bfg^N$ satisfies:
$$ D^\top\bfg^N= \sgn\gtau\,\cdot\vf_{N\,\#}\bigl( \gk_{N-1}\,\gn_{N-2}\,d\calL^1\pri]a,b[\bigr)  $$
where, we recall, $\vf_N(s):=\int_a^s\Vert\dot\gn_N(\l)\Vert\,d\l=\int_a^s|\gtau(\l)|\,d\l $.
%
%
\eex
{\large\sc The milder case.} Assume now that the open curve $\gc$ is mildly smoothly turning at order $j+1$ for some $1\leq j\leq N$, see Definition~\ref{Dms}.
This time, by Proposition~\ref{Cmild} we know that
$ \Pi(\gn_j(s))=[\gn_j](\gc)(\vf_j(s))$ for each $s\in[a,b]$, where the transition function $\vf_j:[a,b]\to[0,L_j]$ is bijective and absolutely continuous.
Moreover, on account of Remark~\ref{Rmild}, the $j$-th normal $\gn_j(s)$ is a function of class $C^1$ in each open interval given by a connected component of $]a,b[\sm\SS$, where $\SS$ is a finite set of points $s_i\in]a,b[$, and
$\gn_j(s_i-)=\pm \gn_j(s_i+)$.
\par Therefore, in this case we can only find a (non continuous) lifting $\bfg^j$ of the function $[\gn_j](\gc)$ such that $\bfg^j(\vf_j(s))=\gn_j(s)$ for each $s\in ]a,b[\sm\SS$.
As a consequence, formula \eqref{Dgj} holds true, but this time equality \eqref{dergj} is satisfied on $]a,b[\sm\SS$,
and it turns out that $\dot\bfg^j:[0,L_j]\to\SN$ is a {\em special function of bounded variation}.
\par More precisely, the distributional derivative of the function $\dot\bfg^j$ decomposes into the absolutely continuous and singular components (w.r.t. the Lebesgue measure $\calL^1$)
$$D\dot\bfg^j=D^a\dot\bfg^j+D^s\dot\bfg^j $$
where, arguing as in formula \eqref{repr}, we have:
$$ \lan D^a\dot\bfg^j,\x\ran=\int_a^{b} {d\over ds}\,{\dot\gn_j(s)\over \Vert\dot\gn_j(s)\Vert}\bullet \x(\vf_j(s))\,d\calL^1(s)
$$
and the singular component is concentrated at the points $s_i\in \SS$, namely:
$$ \lan D^s\dot\bfg^j,\x\ran =\sum_{s_i\in \SS}\Bigl[ {\dot\gn_j(s_i+)\over \Vert\dot\gn_j(s_i+)\Vert}-{\dot\gn_j(s_i-)\over \Vert\dot\gn_j(s_i-)\Vert}\Bigr]\,\x(\vf_j(s_i))\,. $$
%
%
%
%
\par However, by the formulas \eqref{jumps} it turns out that the jumps appearing in the singular component of the
measure derivative $D\dot\bfg^j$, are produced by couples of antipodal point in the Gauss sphere $\SN$.
As a consequence, they cannot be seen in the projective space $\RPN$, and the projected function
$\ds s\mapsto \Pi\circ{\dot\gn_j(s)\over \Vert\dot\gn_j(s)\Vert}$ is continuous in $]a,b[$ and differentiable outside the singular points $s_i\in\SS$.
\par In conclusion, coming back to the weak $j$-th normal $[\gn_j](\gc)=\Pi\circ\bfg^j$, where $\Pi:\SN\to\RPN$ is the canonical projection,
similarly to the smoothly turning case, if the curve $\gc$ is mildly smoothly turning at order $j+1$, then the distributional derivative
of the arc-length derivative of $[\gn_j](\gc)$ is an absolutely continuous measure, and on account of \eqref{meas} we may conclude with the formula:
$$ D {d\over dt}[\gn_j](\gc)=\vf_{j}^\#\wid\m_j\,,\qquad \wid\m_j:={d\over ds}\Bigl(\Pi\circ\,{\dot\gn_j(s)\over \Vert\dot\gn_j(s)\Vert}\Bigr)\,\calL^1\pri]a,b[ $$
that makes sense by means of an isometric embedding of $\RPN$ into some Euclidean space.
\appendix\section{Proof of Proposition~\ref{PGM4}}
Assuming $N=3$, according to the notation from \eqref{vjh}, the fifth order expansions of $\gc$ at $s$ give:
%
$$  \ds \gv_0(h)=\dot\gc+\frac{\gc^{(3)}}6\,h^2+\ga\,h^4+\go(h^4) $$
$$\ds \gv_1(h)=-\dot\gc+2\gc^{(2)}\,h-\frac{13}6\,{\gc^{(3)}}\,h^2+\frac 53\,\gc^{(4)}\,h^3-\gb\,h^4+\go(h^4) $$
$$
\ds \gv_2(h)=\dot\gc+2\gc^{(2)}\,h+\frac{13}6\,{\gc^{(3)}}\,h^2+ \frac 53\,\gc^{(4)}\,h^3+\gb\,h^4+\go(h^4) $$
$$
\ds \gv_3(h)=-\dot\gc+4\gc^{(2)}\,h-\frac{49}6\,{\gc^{(3)}}\,h^2+\frac {34}3\,\gc^{(4)}\,h^3+\go(h^3) $$
where $\ga$ and $\gb$ depend on $\gc^{(5)}(s)$. We thus get:
$$\Vert\gv_0(h)\Vert^2=1-\frac{\Vert\gc^{(2)}\Vert^2}3\,h^2+\Bigl( 2 \ga\bullet\dot\gc+\frac 1{36}\,\Vert\gc^{(3)}\Vert^2\Bigr)\,h^4+o(h^4) $$
and
$$ \Vert\gv_0(h)\Vert^{-2}=1+\frac{\Vert\gc^{(2)}\Vert^2}3\,h^2+\Bigl( \frac 19\,\Vert\gc^{(2)}\Vert^4-\frac 1{36}\,\Vert\gc^{(3)}\Vert^2-2 \ga\bullet\dot\gc\Bigr)\,h^4+o(h^4) $$
whence \eqref{th} holds.
We also have
$$ \gv_1(h)\bullet\gv_0(h)=-1+\frac 73\,\Vert\gc^{(2)}\Vert^2h^2+\frac 13\,\Bigl( \gc^{(3)}\bullet\gc^{(2)} + 5\gc^{(4)}\bullet\dot\gc\Bigr)\,h^3 -a\,h^4 +o(h^4) $$
where
$$ a:= \frac{13}{36}\,\Vert\gc^{(3)}\Vert^2+(\ga+\gb)\bullet\dot\gc $$
and hence
$$ \frac{\gv_1(h)\bullet\gv_0(h)}{\Vert\gv_0(h)\Vert^2}=-1+2\Vert\gc^{(2)}\Vert^2h^2+\frac 13\,\Bigl( \gc^{(3)}\bullet\gc^{(2)} + 5\gc^{(4)}\bullet\dot\gc\Bigr)\,h^3
-b\,h^4+ o(h^4)$$
where
$$ b:= \frac{1}{3}\,\Vert\gc^{(3)}\Vert^2+(\gb-\ga)\bullet\dot\gc-\frac{2}{3}\,\Vert\gc^{(2)}\Vert^4 $$
that gives
$$ \ba{rl} \gN_1(h)= & \ds 2\gc^{(2)}h-2\bigl(\Vert\gc^{(2)}\Vert^2\dot\gc+\gc^{(3)}\bigr)\,h^2+ \frac 13\,\Bigl( 5 \gc^{(4)} - 5(\gc^{(4)}\bullet\dot\gc)\,\dot\gc- (\gc^{(3)}\bullet\gc^{(2)})\dot\gc\Bigr)\,h^3 \\
& \ds +\Bigl( b\,\dot\gc -\frac 13\,\Vert\gc^{(2)}\Vert^2\gc^{(3)}+\ga-\gb\Bigr)\,h^4+\go(h^4)\,. \ea
 $$
As a consequence, we get
$$  \ds \Vert\gN_1(h)\Vert^2=4\Vert\gc^{(2)}\Vert^2h^2- 8 {\gc^{(3)}\bullet\gc^{(2)}}\,h^3+ 4\Bigl(\Vert \gc^{(3)}\Vert^2-\Vert\gc^{(2)}\Vert^4+\frac 53\,
\gc^{(4)}\bullet\gc^{(2)} \Bigr)\,h^4+ o(h^4) $$
whence
$$
\ds \Vert\gN_1(h)\Vert^{-2}=\frac 1{4\Vert\gc^{(2)}\Vert^2h^2}\Bigl[ 1+ 2\frac{\gc^{(3)}\bullet\gc^{(2)}}{\Vert\gc^{(2)}\Vert^2}\,h+
\Bigl( 4\frac{(\gc^{(3)}\bullet\gc^{(2)})^2}{\Vert\gc^{(2)}\Vert^4}
-\frac 53\,\frac{ \gc^{(4)}\bullet\gc^{(2)} }{\Vert\gc^{(2)}\Vert^2}-\frac{\Vert\gc^{(3)}\Vert^2}{\Vert\gc^{(2)}\Vert^2}+ {\Vert\gc^{(2)}\Vert^2}
\Bigr)\,h^2+ o(h^2)\Bigr] $$
and definitely \eqref{N1h} holds, where
\beq\label{d}
{\bf d}:= -\frac 16\,\frac {\gc^{(3)}\bullet\gc^{(2)}}{\Vert\gc^{(2)}\Vert}\,\gt+\O\,\gn_1+\Bigl(\frac 56\,\frac {\gc^{(4)}\bullet\gc^{(3)\perp}}{\Vert\gc^{(2)}\Vert\,\Vert\gc^{(3)\perp}\Vert}- \frac {\gc^{(3)}\bullet\gc^{(2)}}{\Vert\gc^{(2)}\Vert^3}\,\Vert\gc^{(3)\perp} \Vert \Bigr)\,\gn_2+
\frac 56\,\frac {\Vert\gc^{(4)\perp}\Vert}{\Vert\gc^{(2)}\Vert}\,\gn_3
\eeq
with the coefficient $\O$ of $\gn_1$ equal to
\beq\label{Om} \O:= \frac {\bigl(\gc^{(3)}\bullet\gc^{(2)}\bigr)^2}{\Vert\gc^{(2)}\Vert^4}\,\Bigl(\frac 32\,\Vert\gc^{(2)}\Vert^2-1\Bigr)+\frac 12\,\Vert\gc^{(2)}\Vert^2-\frac 12\,\frac{\Vert \gc^{(3)} \Vert^2}{\Vert\gc^{(2)}\Vert^2}\,.
\eeq
%
\par Moreover, in order to compute $\gN_2(h)$, we check:
$$ \gv_2(h)\bullet\gv_0(h)=1-\frac 73\,\Vert\gc^{(2)}\Vert^2h^2+\frac 13\bigl( \gc^{(3)}\bullet\gc^{(2)}+5\gc^{(4)}\bullet\dot\gc\bigr)+o(h^3) $$
$$ \frac{\gv_2(h)\bullet\gv_0(h)}{\Vert\gv_0(h)\Vert^2}=1-2\Vert\gc^{(2)}\Vert^2h^2+\frac 13\bigl( \gc^{(3)}\bullet\gc^{(2)}+5\gc^{(4)}\bullet\dot\gc\bigr)+o(h^3)$$
and hence
$$ -\frac{\gv_2(h)\bullet\gv_0(h)}{\Vert\gv_0(h)\Vert^2}\,\gv_0(h)=
-\dot\gc+\Bigl(2\Vert\gc^{(2)}\Vert^2\,\dot\gc-\frac {1}6\,\gc^{(3)}\Bigr)\,h^2
-\frac 13\,\bigl( \gc^{(3)}\bullet\gc^{(2)}+5\gc^{(4)}\bullet\dot\gc \bigr)\,\dot\gc\,h^3 +\go(h^3)\,. $$
Furthermore,
$$ \gv_2(h)\bullet\gN_1(h)=4\Vert\gc^{(2)}\Vert^2h^2
+\Bigl(\frac{5}{3}\,\gc^{(4)}\bullet\gc^{(2)} + \Vert \gc^{(2)}\Vert^4 - \Vert\gc^{(3)}\Vert^2 \Bigr)\,h^4+o(h^4)
$$
so that
$$\bigl( \gv_2(h)\bullet\gN_1(h)\bigr)\,\gN_1(h)=4\Vert\gc^{(2)}\Vert\,h^2\Bigl\{ 2\gc^{(2)}h-2\bigl(\Vert\gc^{(2)}\Vert^2\dot\gc+\gc^{(3)}\bigr) \,h^2+{\bf A}\,h^3+o(h^3)\Bigr\}  $$
where
$$ {\bf A}:=\frac 13\,\Bigl( 5 \gc^{(4)} - 5(\gc^{(4)}\bullet\dot\gc)\,\dot\gc- (\gc^{(3)}\bullet\gc^{(2)})\dot\gc\Bigr)
+2\,\Bigl(\Vert\gc^{(2)}\Vert^2-\frac{\Vert\gc^{(3)}\Vert^2}{\Vert\gc^{(2)}\Vert^2} \Bigr)\,\gc^{(2)}+\frac{10}3\,\frac{\gc^{(4)}\bullet\gc^{(2)}}{\Vert\gc^{(2)}\Vert^2}\,\gc^{(2)}\,.
$$
We thus obtain:
$$ \frac{ \gv_2(h)\bullet\gN_1(h)}{\Vert\gN_1(h)\Vert^2}\,\gN_1(h) = 2\gc^{(2)}\,h+\Bigl(4\frac{\gc^{(2)}\bullet\gc^{(3)}}{\Vert\gc^{(2)}\Vert^2}\,\gc^{(2)}
-2\bigl(\Vert\gc^{(2)}\Vert^2\dot\gc+\gc^{(3)}\bigr)\Bigr)\,h^2+{\bf B}\,h^3+\go(h^3)\,, $$
where
$$ \ba{l} \ds {\bf B}:=\frac 13\,\Bigl( 5 \gc^{(4)} - 5(\gc^{(4)}\bullet\dot\gc)\,\dot\gc- 13(\gc^{(3)}\bullet\gc^{(2)})\dot\gc\Bigr)
 \\ \ds \qquad
+4\,\Bigl(\Vert\gc^{(2)}\Vert^2-\frac{\Vert\gc^{(3)}\Vert^2}{\Vert\gc^{(2)}\Vert^2} +2\,\frac{\bigl(\gc^{(3)}\bullet\gc^{(2)}\bigr)^2}{\Vert\gc^{(2)}\Vert^4} \Bigr)\,\gc^{(2)}  -4\,\frac{\gc^{(3)}\bullet\gc^{(2)}}{\Vert\gc^{(2)}\Vert^2}\,\gc^{(3)}\,. \ea
$$
Putting the terms together, we get:
$$ \gN_2(h)=4\gc^{(3)\perp} h^2+4\,{\bf D}\,h^3+\go(h^3)  $$
where
$$ {\bf D}:=\bigl(\gc^{(3)}\bullet\gc^{(2)}\bigr)\dot\gc
+\Bigl(\frac{\Vert\gc^{(3)}\Vert^2}{\Vert\gc^{(2)}\Vert^2}-\Vert\gc^{(2)}\Vert^2 \Bigr)\,\gc^{(2)} +\frac{\gc^{(3)}\bullet\gc^{(2)}}{\Vert\gc^{(2)}\Vert^2} \,\gc^{(3)} -2\,\frac{\bigl(\gc^{(3)}\bullet\gc^{(2)}\bigr)^2}{\Vert\gc^{(2)}\Vert^4}\,\gc^{(2)}
$$
and definitely \eqref{N2h} holds, where in terms of the orthonormal basis $(\gt,\gn_1,\gn_2,\gn_3)$ we obtain the formula \eqref{D} for ${\bf D}$.
\par Finally, formula \eqref{N3h} follows by arguing as in the proof of Proposition~\ref{PGMN}. In fact, the Gram-Schmidt procedure yields that $(\gt(h),\gn_1(h),\gn_2(h),\gn_3(h))$ is an orthonormal basis of $\gR^4$, whence $\gn_3(h)=\gn_3+\go(1)$.
%
\par
More precisely, we have $\gn_3(h)=\pm \ast(\gt(h)\wdg\gn_1(h)\wdg\gn_2(h))$, where $\ast$ is the Hodge operator in $\gR^4$, whereas
$\ast(\gt\wdg\gn_1\wdg\gn_2)=\pm\gn_3$, with the same sign $\pm$ in the previous two formulas, by our choice in \eqref{vjh}.
Using that
$$ \gt(h)=\gt+\go(h)\,,\quad \gn_1(h)=\gn_1+\a\,\gn_2\,h+\go(h)\,,
%
%
\quad \gn_2(h)=\gn_2+(\be\,\gt+\g\,\gn_1)\,h+\go(h)
$$
for some real numbers $\a,\be,\g\in\gR$, we get $\gt(h)\wdg\gn_1(h)=\gt\wdg\gn_1+\a\,\gt\wdg\gn_2\,h+\go(1)\wdg\go(1)\,h$ and hence
$\gt(h)\wdg\gn_1(h)\wdg\gn_2(h)=\gt\wdg\gn_1\wdg\gn_2+\go(1)\wdg\go(1)\wdg\go(1)\,h$, whence actually \eqref{N3h} holds true, as required.
\adl\par\noindent
{\bf Acknowledgements.}
The research of D.M. was partially supported by the GNAMPA of INDAM.
The research of A.S. was partially supported by the GNSAGA of INDAM.
 \section*{Conflict of interest}
 The authors declare that they have no conflict of interest.
\end{document}